\documentstyle[11pt]{article}

\input amssym.def \input amssym
\font\tencmmib=cmmib10 \font\sevencmmib=cmmib7 \font\fivecmmib=cmmib5
\newfam\cmmibfam \textfont\cmmibfam=\tencmmib
\scriptfont\cmmibfam=\sevencmmib \scriptscriptfont\cmmibfam=\fivecmmib

 \font\teneufm=eufm10 at 12 pt
\font\seveneufm=eufm9 \font\fiveeufm=eufm7
\newfam\eufmfam \textfont\eufmfam=\teneufm
\scriptfont\eufmfam=\seveneufm \scriptscriptfont\eufmfam=\fiveeufm
\def\frak#1{{\fam\eufmfam\relax#1}}

\font\script=eusm10 at 12 pt

\font\helv=cmssbx10 at 12 pt
\def\isom{\ifmmode\ \cong\ \else$\isom$\fi}
\def\scsi{\scriptsize}
\def\scst{\scriptstyle}
\def\scstst{\scriptscriptstyle}
\def\aa{\frak a}  \def\pp{\frak
    p}  \def\Gg{\frak g}
 \def\hh{\frak h} \def\mm{\frak
    m} \def\nn{\frak n} 
\def\tt{\frak t} \def\kk{\frak k} \def\ll{\frak
    l} \def\Ss{\frak s} 
\def\uu{\frak u}   \def\zz{\frak z}

\def\AA{\frak A} \def\BB{\frak B}   
 \def\JJ{\frak J}

   \def\CC{\frak C}

 \def\fC{{\Bbb C}}
 
 \def\fQ{{\Bbb Q}}
\def\fR{{\Bbb R}} \def\fH{{\Bbb H}}

\def\bfb{\hbox{\scsi \bf b}}
\def\bfn{\hbox{\scsi \bf n}}
\def\hbf#1{\hbox{\scsi\bf #1}}

\def\bs{\backslash}

\def\gD{\Delta} \def\gd{\delta} \def\gG{\ifmmode {\Gamma} \else$\gG$\fi}
 \def\gs{\sigma} \def\gS{\Sigma} \def\gz{\zeta}
 \def\gt{\theta}  \def\ga{\alpha}
\def\gb{\beta} \def\gk{\kappa}  
    \def\gr{\rho}
 \def\ge{\varepsilon}  

 \def\-#1{\overline{#1}}
\def\~#1{\tilde{#1}}    \def\cO{\ifmmode {\cal O} \else$\cO$\fi}

\def\cD{{\cal D} }

\def\cE{{\cal E} }

\def\cK{{\cal K} }
\def\cL{{\cal L} }

\def\cN{{\cal N} }
\def\cO{{\cal O} }

\def\cR{{\cal R} }

\def\cZ{{\cal Z} }

\def\scD{\hbox{{\script D}} }

\def\scH{\hbox{{\script H}} }

\def\scJ{\hbox{{\script J}} }

\def\scL{\hbox{{\script L}} }

\def\scS{\hbox{{\script S}} }

 \def\ad{\hbox{ad}} \def\dim{\hbox{dim}}

  \def\Ker{\hbox{Ker}}
 \def\rank{\hbox{rank}} \def\span{\hbox{span}}

\def\isom{\ifmmode\ \cong\ \else$\isom$\fi}

\def\nni{\supset}  \def\ovv{\ifmmode
  \overline{\gG \bs \cal D}\ \else$\ovv$\fi}

\def\xg{\ifmmode {X_{\gG}} \else$\xg$\fi} \def\xgeq{\ifmmode {\xg =
    \gG\bs\cD} \else$\xgeq$\fi}
\def\xgs{\ifmmode {X_{\gG}^*} \else$\xgs$\fi}
\def\xgc{\ifmmode {\overline{X}_{\gG}} \else$\xgc$\fi}

\def\cIf{\ifmmode {\cal I}_5%
  \else$\cIf$\/\fi} \def\cIt{\ifmmode {\cal I}_{10} \else$\cIt$\fi}
\def\JG{\ifmmode {G_{25,920}} \else$\JG$\fi}

\def\pts{\ifmmode \hbox{\helv P}^2_6 \else $\pts$\fi} \def\pthrees{\ifmmode
  \hbox{\helv P}^3_6 \else $\pthrees$\fi}

\def\pos{\ifmmode \hbox{\helv P}^1_6 \else $\pos$\fi}

\def\hra{\hookrightarrow} 
\def\lra{\longrightarrow} 
 
 \def\ra{\rightarrow} \def\Ra{\Rightarrow}
 \def\bs{\ifmmode {\setminus} \else$\bs$\fi}

\def\inn{\subset} \def\nni{\supset}

\def\ende{\hfill $\Box$ \vskip0.25cm }

\newtheorem{theorem}{Theorem}[section]

\newtheorem{Lemma}{Lemma}[theorem]

\newtheorem{lemma}[theorem]{Lemma}
\newtheorem{proposition}[theorem]{Proposition}
\newtheorem{corollary}[theorem]{Corollary}

\newenvironment{definition}{\refstepcounter{theorem}
\vskip.2cm\noindent{\bf Definition \thesection.\arabic{theorem}
  }}{\vskip.2cm\noindent}

\newenvironment{examples}{\refstepcounter{theorem}
\vskip.2cm\noindent{\bf Examples \thesection.\arabic{theorem}:
  }}{\vskip.2cm\noindent}

\begin{document}

\title{Symmetric subgroups of rational groups of hermitian type}
\author{Bruce Hunt \\ FB Mathematik, Universit\"at \\ Postfach 3049 \\
  67653 Kaiserslautern}
\maketitle

\begin{center}
{\LARGE\bf Introduction}
\end{center}

A rational group of hermitian type is a $\fQ$-simple algebraic group $G$
such that the symmetric space $\cD$ of maximal compact subgroups of the
real Lie group $G(\fR)$ is a hermitian symmetric space of the non-compact
type. There are two major
classes of subgroups of $G$ of importance to the geometry of $\cD$ and to
arithmetic quotients $X_{\gG}=\gG\bs \cD$ of $\cD$: {\it parabolic}
subgroups and {\it reductive} subgroups. The former are connected
with the boundary of the domain $\cD$, while the latter are connected with
submanifolds $\cD'$ of the {\it interior} of $\cD$, {\it a priori} just
symmetric
spaces. Under certain assumptions on the reductive subgroup,
$\cD'\inn \cD$ is a {\it
  holomorphic symmetric embedding}, displaying $\cD'$ as a sub-hermitian
symmetric space. It is these latter groups we study in this
paper, and we call them quite generally ``symmetric subgroups''; these are
the subgroups occuring in the title.

The purpose of this paper is to prove an existence result of the following
kind: given a maximal $\fQ$-parabolic $P\inn G$, which is the stabilizer of
a boundary component $F$, meaning that $P(\fR)=N(F)$,
there exists a $\fQ$-subgroup
$N\inn G$, which is a symmetric subgroup defining a subdomain
 $\cD_N\inn \cD$ such
that: $F$ is a boundary component of $\cD_N$. To formulate this precisely,
we make the following assumptions on $G$ which are to be in effect
throughout the paper: $G$ is isotropic, and $G(\fR)$ is not a product of
groups of type $SL_2(\fR)$.
The basis of our formulation is the notion of {\it incident} parabolic and
symmetric subgroups. We first define this for the real groups. Fix a
maximal $\fR$-parabolic $P\inn G(\fR)$, and let as above $F\inn \cD^*$
denote the corresponding boundary component stabilized by $P$. First assume
that $F$ is positive-dimensional. Let $N\inn G(\fR)$ be a symmetric
subgroup; we shall say $N$ and $P$ are {\it incident}, if the following
conditions are satisfied.
\begin{itemize}\item[1)] $N$ has maximal $\fR$-rank, that is,
  $\rank_{\fR}N=\rank_{\fR}G$.
\item[2)] $N$ is a maximal symmetric subgroup.
\item[3)] $N=N_1\times N_2$, where $N_1\inn P$ is a hermitian Levi factor of
  $P$ for some Levi decomposition (these terms are explained in more detail
  in the text).
\end{itemize}
{}From work of Satake and Ihara one knows, up to conjugation, all subgroups
$N$ fulfilling the above; these results are collected in Table \ref{T1},
and such a group is uniquely determined up to conjugation (by an element of
$K\inn G(\fR)$, a maximal compact subgroup). It is trickier to define
something similar for the case $\dim(F)=0$.
This is easiest to see by considering the geometric equivalent of the
conditions above. They state that the subdomain $\cD_N$ defined by $N$ has
the same $\fR$-rank as $\cD$, is a maximal subdomain, and finally, that
$\cD_N$ contains the given boundary component $F$ in {\it its} boundary,
$F\inn \overline{\cD}_N$, and furthermore that $\cD_N=\cD_1\times \cD_2$,
where $\cD_1$ is (as an abstract hermitian symmetric space) isomorphic to
$F$. The condition 1) is still a natural one to assume, but 2) is in
general too strong, while 3), in the case of $\dim(F)=0$, seems to imply
$\cD_N$ should be irreducible and contain $F$ as a boundary component. We
then replace, for $\dim(F)=0$, the conditions 2) and 3) by 2') and 3') below:

\begin{itemize}\item[2')] $N$ is a maximal subgroup of tube type, i.e.,
  such that $\cD_N$ is a tube domain.
\item[3')] $\cD_N$ is maximal irreducible and contains $F$ as a
boundary component.
\end{itemize}
The subgroups $N$ fulfilling 1), 2), 3') or 1), 2'), 3') are also known and
are listed in Table \ref{T2}. But there simply do not always exist such
subgroups; more precisely, for the domains

$(\cE\cD)\hspace*{5cm} \bf I_{\hbox{\scsi\bf q,q}},\ II_{\hbox{\scsi\bf
    n}},\ n \hbox{ even},\ III_{\hbox{\scsi\bf n}},$

\noindent there are no subgroups $N$ fulfilling the above.
So in these cases, we introduce the conditions
\begin{itemize}\item[2'')] $N$ is {\it minimal}, subject to 1).
\item[3'')] $\cD_N$ contains $F$ as a boundary component.
 \end{itemize}
It is clear that condition 1) and 2'') together imply that the domain
$\cD_N$ is a polydisc, $\cD_N\cong (\cD_{1,1})^t$, where $t=\rank_{\fR}G$,
and $\cD_{1,1}=SL_2(\fR)/K$ is the one-dimensional disc.
Altogether these conditions insure that for any maximal
$\fR$-parabolic $P\inn G(\fR)$, there is a finite set of incident symmetric
subgroups $N\inn G(\fR)$, such that any subgroup incident to $P$ is
isomorphic to one of them. Moreover, if $F$ is positive-dimensional, $N$ is
unique up to isomorphism.

Recall now that
given $G(\fR)$ of hermitian type, $\rank_{\fR}G(\fR)=t$, one can find $t$
strongly orthogonal roots $\mu_1,\ldots, \mu_t$ which determine a unique
maximal $\fR$-split torus $A\inn G(\fR)$, the corresponding root system
$\Phi(A,G)$, as well as a set of simple
$\fR$-roots. The closed, symmetric set of roots
$\Psi=\{\pm\{\mu_1\},\ldots,\pm\{\mu_t\}\}$ then determines a (unique)
$\fR$-subgroup $N_{\Psi}\inn G(\fR)$, which is isomorphic to
$(SL_2(\fR))^t$. Conversely, the maximal $\fR$-split torus $A$ defines the
root system $\Phi(A,G)$, and in this root system there is a good unique
choice for the strongly orthogonal roots $\mu_1,\ldots,\mu_t$ (see section
\ref{section1.1}), so $A$ determines $\Psi$ and $\Psi$ determines $A$.
For each of the symmetric groups above, it is natural to
consider also the condition
\begin{itemize}\item[4)] Assume the parabolic $P$ is standard with respect
  to $A$. Then $N_{\Psi}\inn N$, where $\Psi$ and $A$ determine one another
  in the manner just described.
\end{itemize}
We then pose the following problem. Given the rational group of hermitian
type $G$, and given a maximal $\fQ$-parabolic $P\inn G$, find a symmetric
$\fQ$-subgroup $N\inn G$, such that $P$ and $N$ are {\it incident}, meaning
that $N(\fR)$ and $P(\fR)$ are incident in $G(\fR)$ is the above sense,
i.e., $N(\fR)$ fulfills 1), 2) and 3) (if $P(\fR)=N(F)$, $\dim(F)>0$ ), 1),
2') and 3') or 1), 2'') and 3'') in the corresponding cases $\dim(F)=0,
\cD\not\in (\cE\cD)$ respectively $\dim(F)=0, \cD\in (\cE\cD)$.

The main result of the paper is that the problem just posed can be solved
in almost all cases, specified in the following theorem.

\vspace*{.2cm}
\noindent{\bf Main Theorem} {\it Let $G$ be $\fQ$-simple of hermitian type
  subject to the restrictions above ($G(\fR)$ not a product of
  $SL_2(\fR)$'s and $G$ isotropic), $P\inn G$ a maximal
  $\fQ$-parabolic. Then there exists a reductive $\fQ$-subgroup $N\inn G$
  such that $(P,N)$ are incident, in the sense just defined, with the
  following exceptions:

In the case of zero-dimensional boundary components,
\begin{itemize}
\item Index $C^{(2)}_{2n,n}$, that is, $G$ is isogenous to the group
  $Res_{k|\fQ}G'$, $G'=SU(V,h)$, the special unitary group for a hermitian form
  $h:V\times V\lra D$, where $V$ is a $2n$-dimensional right $D$-vector
  space, $D$ a totally indefinite quaternion division algebra which is central
  simple over $k$ ($k$ a totally real number field), and the Witt index
  of $h$ is $n$.
\end{itemize}
In this case the $k$-subgroups fulfilling 1), 2'') (over $k$) and 3'') and
the corresponding domains are
\begin{itemize}\item $n>1$: $N=Res_{k|\fQ}N'$, $N'\cong SU(V_1,h_{|V_1})\times
  \cdots \times
  SU(V_n, h_{|V_n})$, where $V_i$ is a hyperbolic plane, $V\cong V_1\oplus
  \cdots\oplus V_n$, and then $N'(\fR)\cong Sp(4,\fR)\times \cdots \times
  Sp(4,\fR)$, and the domain $\cD_{N'}$ is of type $\bf III_{\hbox{\scsi\bf
      2}}\times \cdots \times III_{\hbox{\scsi\bf 2}}$.
\item $n=1$: $N=Res_{k|\fQ}N'$, $N'\cong SU(V_L,h_{|V_L})$, where
  $V_L\inn V$ is a plane defined
  by the inclusion of an imaginary quadratic extension $L$ of $k$, $L\inn D$
  (viewing $V$ as a $k$-vector space), and then $N'(\fR)\cong SL_2(\fR)$ and
  the domain $\cD_{N'}$ is of type $\bf III_{\hbox{\scsi\bf 1}}$, a disc
  (in particular {\em not} fulfilling 1)).
\end{itemize}
Moreover, there is a $\fQ$-subgroup $N$ fulfilling also condition 4) with
the exception of the index $C^{(2)}_{2,1}$.}

\vspace*{.2cm}
Hence, for indices $C^{(2)}_{2n,n}$ with $n>1$ we can find $N$
fulfilling 1) and 3''), but neither 2')
nor 2''), but rather only 2'') over $k$, that is, requiring $N$ to be
minimal {\it over k}. In the case $n=1$, it would seem that there
are no subgroups
satisfying condition 1) at all; at any rate we could find none. This case
was considered in detail in \cite{hyp}. As a final remark, we could
eliminate the exceptional status of the cases $n>1$ by changing 2'')
accordingly, but we feel these groups are indeed exceptions.

The main application of the above theorem is to arithmetic quotients of
bounded symmetric domains. Let $\gG\inn G(\fQ)$ be an arithmetic subgroup,
$X_{\gG}=\gG\bs \cD$ the arithmetic quotient. Since the subgroups $N\inn G$
of the theorem are defined over $\fQ$, $\gG_N:=\gG\cap N(\fQ)$ is an
arithmetic subgroup in $N$, and $X_{\gG_N}:=\gG_N\bs \cD_N$ is an
arithmetic quotient. Clearly $X_{\gG_N}\inn X_{\gG}$, and, viewing
$X_{\gG}$ as an algebraic variety, $X_{\gG_N}$ is a subvariety, which we
call a {\it modular} subvariety. More precisely, one has the Baily-Borel
embedding of the Satake compactification $X_{\gG}^*$, which makes
$X_{\gG}^*$ an algebraic variety and $X_{\gG}$ a Zariski open subset of a
projective variety. By a theorem of Satake \cite{S2}, the natural
inclusion $X_{\gG_N}\inn X_{\gG}$ extends to an inclusion of the
Baily-Borel compactifications, displaying $X_{\gG_N}^*$ as a (singular)
algebraic subvariety of the algebraic variety $X_{\gG}^*$. The notion of
incidence then manifests itself in the following way. A boundary component
$V_i\inn X_{\gG}^*$ (here $V_i$ is itself an arithmetic quotient) and
modular subvariety $X_{\gG_N}^*$ are incident, if $V_i$ is a boundary
component of $X_{\gG_N}^*$. These matters will be taken up elsewhere.

Let us now sketch the idea of proof of the Main Theorem. We consider the
different possibilities for the $\fQ$-rank of $G$. By assumption, $G$ is
isotropic, so $\rank_{\fQ}G\geq 1$. We split the possible cases into the
following three items, corresponding roughly to increasing order of
difficulty:
\begin{itemize}\item[1)] $G$ is split over $\fR$ (explained below).
\item[2)] $G$ is not split over $\fR$, but $\rank_{\fQ}G\geq 2$.
\item[3)] $\rank_{\fQ}G=1$.
\end{itemize}
The first case occurs as follows. Since $G$ is $\fQ$-simple, there is a
totally real number field $k$ and an {\it absolutely} simple $k$-group $G'$
such that $Res_{k|\fQ}G'=G$. Let $S$ be a maximal $k$-split torus of $G'$,
$A\inn G'(\fR)$ a maximal $\fR$-split torus, which one may assume contains
$S$, $S\inn  A$. Then $G'$ is {\it split over $\fR$}, if $S=A$, i.e., a
maximal $k$-split torus is also a maximal $\fR$-split torus. In this case
the proof of the main theorem is more or less trivial, and is carried out
in \S4. Case 2) is already more subtle; we combine constructions based on
the correspondence between semisimple Lie algebras (of the classical type)
and central simple algebras with involution, as presented in \cite{W}, with
the known classification of $k$-indices of absolutely simple groups, as
presented in \cite{T}. This is carried out in \S5, and in this as well as
the next section most of the work is necessary for the case of
zero-dimensional boundary components. The most interesting
case, done in the last paragraph, is the case of rank one. Here it turns
out that similar arguments as above apply in the cases that the boundary
component $F$ is positive-dimensional; but in the case of zero-dimensional
boundary components, new arguments are required. In particular, the case of
{\it hyperbolic planes} is fundamental, and this was studied in detail in
\cite{hyp}. Drawing on the results of that paper we can complete the proof
of the theorem.

\section{Real parabolics of hermitian type}
\subsection{Notations}\label{section1.1}
In this paper we will basically adhere to the notations of \cite{BB}. In
the first two paragraphs $G$ will denote a real Lie group; later $G$ will
be a $\fQ$-group of hermitian type. We assume $G$ is reductive, connected
and with compact center; $K\inn G$ will denote a maximal compact subgroup,
$\cD=G/K$ the corresponding symmetric space. Throughout this paper we will
assume $G$ is of {\it hermitian type}, meaning that $\cD$ is a hermitian
symmetric space, hence a product $\cD=\cD_1\times \cdots \times \cD_d$ of
irreducible factors, each of which we assume is {\it non-compact}. Let
$\Gg=\kk+\pp$ denote a Cartan decomposition of the Lie algebra of $G$,
$\Gg_{\fC}=\kk_{\fC}+\pp^++\pp^-$ the decomposition of the complexified Lie
algebra, with $\pp^{\pm}$ abelian subalgebras (and $\pp_{\fC}=\pp^+\oplus
\pp^-$). Chooosing a Cartan subalgebra $\hh\inn \Gg$, the set of roots of
$\Gg_{\fC}$ with respect to $\hh_{\fC}$ is denoted
$\Phi=\Phi(\hh_{\fC},\Gg_{\fC})$. As usual, we choose root vectors
$E_{\ga}\in \Gg^{\ga}$ such that the relations
$$[E_{\ga},E_{-\ga}]=H_{\ga}\in \hh_{\fC},\quad
\ga(H_{\gb})=2{<\ga,\gb>\over <\gb,\gb>},\ \ \ga,\gb\in \Phi,$$
hold.
Complex conjugation maps $\pp^+$ to $\pp^-$, in fact permuting $E_{\ga}$
and $E_{-\ga}$ for $E_{\ga}\in \pp^{\pm}$. Moreover, if
  $\gS^{\pm}:=\{\ga|E_{\ga} \in \pp^{\pm}\}$, then
$$\pp^{\pm}=\span_{\fC}(E_{\ga}),\ \ga\in \gS^{\pm};\quad
\pp=\span_{\fR}(X_{\ga},Y_{\ga}),\ \ga\in \gS^+,$$
where $X_{\ga}=E_{\ga}+E_{-\ga},\ Y_{\ga}=i(E_{\ga}-E_{-\ga})$ (twice the
real and the
(negative of the) imaginary parts, respectively). Let
$\mu_1,\ldots,\mu_t$ denote a maximal
set of strongly orthogonal roots, determined
as in \cite{Helg}: $\mu_1$ is the smallest root in $\gS^+$, and $\mu_j$ is
the smallest root in $\gS^+$ which is strongly orthogonal to
$\mu_1,\ldots,\mu_{j-1}$. This set will be fixed once and for all.

Once this set of strongly orthogonal roots has been chosen, a maximal
$\fR$-split torus $A$ is uniquely determined by
$Lie(A)=\aa=\span_{\fR}(X_{\mu_1},\ldots,X_{\mu_t})$. Then
$\Phi_{\fR}=\Phi(\aa,\Gg)$ will denote the set of $\fR$-roots, and $\Gg$
has a decomposition
$$\Gg=\zz(\aa)\oplus \sum_{\eta\in\Phi_{\fR}}\Gg^{\eta},$$
where $\Gg^{\eta}=\{x\in \Gg | \ad(s)x=\eta(s)x, \forall_{s\in A}\}$. For
each irreducible component of $G$, the set of $\fR$-roots is either of type
$\bf C_{\hbox{\scsi\bf t}}$ or $\bf BC_{\hbox{\scsi\bf t}}$, and of type
$\bf C_{\hbox{\scsi\bf t}}$ $\iff$ the corresponding
domain is a tube domain. If $\xi_i$ denote coordinates on $\aa$ dual to
$X_{\mu_i}$, assuming for the moment $\cD$ to be irreducible, the
$\fR$-roots are explicitly
\begin{equation}\label{e22b.2}
\begin{array}{lcr}\Phi_{\fR}: & \pm(\xi_i\pm \xi_j),\ \pm 2 \xi_i \
(1\leq i\leq t,\ i<j) &  \mbox{(Type $\bf C_{\hbox{\scsi\bf t}}$)} \\
& \pm(\xi_i\pm \xi_j),\ \pm 2 \xi_i, \ \pm\xi_i \
(1\leq i\leq t,\ i<j)
 & \mbox{(Type $\bf BC_{\hbox{\scsi\bf t}}$)} \\
\gD_{\fR}: &  \eta_i=\xi_i-\xi_{i+1}, i=1,\ldots,t-1, \mbox{ and }
\eta_t=2\xi_t \mbox{ (Type $\bf C_{\hbox{\scsi\bf t}}$) },
\eta_t=\xi_t & \mbox{ (Type $\bf BC_{\hbox{\scsi\bf t}}$).}
\end{array}
\end{equation}
Here the simple roots $\gD_{\fR}$ are with respect to the lexicographical
order on the $\xi_i$. A general $\fR$-root system is a disjoint union of
simple root systems. The choice of maximal set of strongly orthogonal roots
determines an order on $\aa$ (the lexicographical order), which
determines, on each simple $\fR$-root system, an order as above; this is
called the {\it canonical order}.

\subsection{Real parabolics}
The maximal $\fR$-split abelian subalgebra $\aa$,
together with the order on it
(induced by the choice of strongly orthogonal roots), determines a unique
nilpotent Lie algebra of $\Gg$, $\nn=\sum_{\eta\in
  \Phi_{\fR}^+}\Gg^{\eta}$. Set $A=\exp(\aa),\ N=\exp(\nn)$, and
\begin{equation}\label{e2.1} B:=\cZ(A)\rtimes N;
\end{equation}
this is a minimal $\fR$-parabolic, the {\it standard} one, uniquely
determined by the choice of strongly orthogonal roots. Every minimal
$\fR$-parabolic of $G$ is conjugate to $B$. Note that, setting
$M=\cZ(A)\cap K$, we have $\cZ(A)=M\times A$, and the group $M$ is the {\it
  semisimple anisotropic kernel} of $G$.

Assume again for the moment that $\cD$ is irreducible, and let $\eta_i,\
i=1,\ldots,t$ denote the simple $\fR$-roots. Set $\aa_b:=\cap_{j\neq b}
\Ker\eta_j,\ b=1,\ldots, t$, a one-dimensional subspace of $\aa$, and
$A_b:=\exp(\aa_b)$, a one-dimensional $\fR$-split subtorus of
$A$. Equivalently, $A_b=\left(\cap_{j\neq b}\Ker\eta_j\right)^0$, where
$\eta_j$ is viewed as a character of $A$. The {\it standard maximal
  $\fR$-parabolic}, $P_b,\ b=1,\ldots, t$, is the group generated by
$\cZ(A_b)$ and $N$; equivalently it is the semidirect product (Levi
decomposition)
\begin{equation}\label{e2.2} P_b=\cZ(A_b)\rtimes U_b,
\end{equation}
where $U_b$ denotes the unipotent radical. The Lie algebra $\uu_b$ of $U_b$
is the direct sum of the $\Gg^{\eta},\ \eta\in \Phi_{\fR}^+,\
\eta_{|\aa_b}\not\equiv 0$. The Lie algebra $\zz(\aa_b)$ of $\cZ(A_b)$ has
a decomposition:
\begin{equation}\label{e2.3} \zz(\aa_b)=\mm_b\oplus \ll_b\oplus
  \ll_b'\oplus \aa_b,\quad \ll_b=\sum_{\eta\in
    [\eta_{b+1},\ldots,\eta_t]}\Gg^{\eta}+[\Gg^{\eta},\Gg^{-\eta}],\ \
  \ll_b'=\sum_{\eta\in
    [\eta_{1},\ldots,\eta_{b-1}]}\Gg^{\eta}+[\Gg^{\eta},\Gg^{-\eta}],
\end{equation}
and $\mm_b$ is an ideal of $\mm$, the Lie algebra of the (semisimple)
anisotropic kernel $M$. Both $\ll_b$ and $\ll_b'$ are simple, and the root
system $[\eta_{b+1},\ldots,\eta_t]$ is of type ${\bf C}_{\hbox{\scsi\bf t-b}}$
or ${\bf
BC}_{\hbox{\scsi\bf t-b}}$, while the root system $[\eta_1,\ldots, \eta_{b-1}]$
is of type
${\bf A}_{\hbox{\scsi\bf b-1}}$. Let $L_b,\ L_b'$ denote the analytic groups
with Lie
algebras $\ll_b$ and $\ll_b'$, respectively, and let $\cR_b=L_b'A_b$, a
reductive group (of type ${\bf A}_{\hbox{\scsi\bf b-1}}$). We call $L_b$ the
{\it hermitian
  factor} of the Levi component and $\cR_b$ the {\it reductive factor}. It
is well known that $L_b$ defines the hermitian symmetric
space which is the $b^{th}$ (standard) boundary component of $\cD$. Indeed,
letting $K_b\inn L_b$ denote a maximal compact subgroup, $\cD_b=L_b/K_b$ is
hermitian symmetric, and naturally contained in $\cD$ as a subdomain,
$\cD_b\inn \cD$. Let $\gz:\cD\lra \pp^+$ be the Harish-Chandra embedding,
and let ${\bf D}=\gz(\cD),\ {\bf D}_b=\gz(\cD_b)$ denote the images; ${\bf
  D}_b$ is a bounded symmetric domain contained in a linear subspace (which
can be identified with $\pp_b^+=\ll_{b,\fC}\cap \pp^+$). Let
$o_b=-(E_{\mu_1}+\cdots+E_{\mu_b}),\ 1\leq b\leq t$;
as the elements $o_b$ are in $\pp^+$, one can consider
the orbits $o_b\cdot G$ and $o_b\cdot L_b$. Since for $g\in L_b$ the action
is described by
$o_bg = o_b +\gz(g)$, one has $o_b\cdot L_b = o_b+\gz(\cD_b)=o_b+{\bf
  D}_b$, and this is the domain ${\bf D}_b$, translated into an affine
subspace ($o_b+\pp_b^+$) of $\pp^+$. One denotes this domain by $F_b:=o_b\cdot
L_b$, and this is the $b^{th}$ {\it standard boundary component} of ${\bf
  D}$. $G$ acts by translations on the various $F_b$, and the images are
the {\it boundary components} of ${\bf D}$; one has
$$\overline{\bf D}={\bf D}\cup \{\hbox{ boundary components
  }\}={\bf D} \cup \left(\cup_{b=1}^to_b\cdot G\right),$$
and $\overline{\bf D}\inn \pp^+$ is the compactification of ${\bf D}$ in
the Euclidean topology. For any boundary component $F$ one denotes by
$N(F),\ Z(F)$ and $G(F)$ the normalizer, centralizer and automorphism group
$G(F)=N(F)/Z(F)$, respectively. Then, letting $U(F)$ denote the unipotent
radical of $N(F)$,
\begin{equation}\label{e3.1} N(F_b)=P_b,\quad U(F_b)=U_b,\quad
  Z(F_b)=Z_b,\quad G(F_b)=L_b,
\end{equation}
where $Z_b$ is a closed normal subgroup of $P_b$ containing every normal
subgroup of $P_b$ with Lie algebra $\zz_b=\mm_b\oplus\ll_b'\oplus
\aa_b\oplus \uu_b$, which is an ideal in $\pp_b$.

Now consider the general case, $\cD=\cD_1\times \cdots \times \cD_d$,
$\cD_i$ irreducible. For each $\cD_i$ we have $\fR$-roots $\Phi_{i,\fR}$,
of $\fR$-ranks $t_i$ and simple $\fR$-roots
$\{\eta_{i,1},\ldots,\eta_{i,t_i}\},\ i=1,\ldots, d$. For each factor we
have standard parabolics $P_{i,b_i}$ $(1\leq b_i\leq t_i$) and standard
boundary components $F_{i,b_i}$. The standard parabolics of $G$ and
boundary components of $\cD$ are then
products
\begin{equation}\label{e3.2} P_{\bfb}=P_{1,b_1}\times \cdots \times
  P_{d,b_d},\quad F_{\bfb}=F_{1,b_1}\times \cdots \times
  F_{d,b_d},\quad ({\bf b}=(b_1,\ldots, b_d)),
\end{equation}
and as above $P_{\bfb}=N(F_{\bfb})$. Furthermore,
\begin{equation}\label{e3.3} G(F_{\bfb})=:L_{\bfb}=L_{1,b_1}\times
  \cdots\times L_{d,b_d}.
\end{equation}
As far as the domains are concerned, any of the boundary components
$F_{i,b_i}$ may be the {\it improper} boundary component $\cD_i$,
which is indicated by setting $b_i=0$. Consequently, $P_{i,0}=L_{i,0}=G_i$
and in (\ref{e3.2}) and (\ref{e3.3})
any ${\bf b}=(b_1,\ldots,b_d),\ 0\leq b_i\leq t_i$ are admissible.

\subsection{Fine structure of parabolics}
For real parabolics of hermitian type one has a very useful refinement of
(\ref{e2.2}). This is explained in detail in \cite{SC} and especially in
\cite{S}, \S III.3-4. First we have the decomposition of $\cZ(A_b)$ as
described above,
\begin{equation}\label{e3.4} \cZ(A_b)=M_b\cdot L_b \cdot \cR_b,
\end{equation}
where $M_b$ is compact, $L_b$ is the hermitian Levi factor, $\cR_b$ is
reductive (of type $\bf A_{\hbox{\scsi \bf b-1}}$), and the product is almost
direct (i.e., the factors have finite intersection). Secondly, the
unipotent radical decomposes,
\begin{equation}\label{e3.5} U_b=\cZ_b\cdot V_b,
\end{equation}
which is a direct product, $\cZ_b$ being the center of $U_b$. The action of
$\cZ(A_b)$ on $U_b$ can be explicitly described, and is the basis for the
compactification theory of \cite{SC}. Before we recall this, let us note
the notations used in \cite{SC} and \cite{S} for the decomposition. In
\cite{SC}, we find
\begin{equation}\label{e4.1} P(F)=(M(F)G_h(F)G_{\ell}(F))\rtimes U(F)\cdot
  V(F),
\end{equation}
and in \cite{S}, where the author
 uses Hermann homomorphisms $\gk:\Ss\ll_2(\fR)\lra \Gg$ to
index the boundary component,
\begin{equation}\label{e4.2} B_{\gk}=\left(G_{\gk}^{(1)}\cdot
    G_{\gk}^{(2)}\right) \rtimes U_{\gk}V_{\gk}.
\end{equation}
In (\ref{e4.1}), $M=M_b,\ G_h=L_b,\ G_{\ell}=\cR_b$, while in (\ref{e4.2}),
$G_{\gk}^{(1)}=M_b\cdot L_b,\ G_{\gk}^{(2)}=\cR_b$ in our notations. The
action can be described as follows (\cite{S}, III \S3-4).
\begin{theorem}\label{t4.1} In the decomposition of the standard parabolic
  $P_b$ (see (\ref{e3.4}) and (\ref{e3.5}))
$$P_b=(M_b\cdot L_b \cdot \cR_b)\rtimes \cZ_b\cdot V_b,$$
the following statements hold.
\begin{itemize}\item[(i)] The action of $M_b\cdot L_b$ is trivial on
  $\cZ_b$, while on $V_b$ it is by means of a symplectic representation
  $\gr:M_b\cdot L_b \lra Sp(V_b,J_b)$, for a symplectic form $J_b$ on
  $V_b$.
\item[(ii)] $\cR_b$ acts transitively on $\cZ_b$ and defines a homogenous
  self-dual (with respect to a bilinear form) cone $C_b\inn \cZ_b$, while
  on $V_b$ it acts by means of a representation $\gs:\cR_b\lra GL(V_b,I_b)$
  for some complex structure $I_b$ on $V_b$.
\end{itemize}
Furthermore the representations $\gr$ and $\gs$ are compatible in a natural
sense.
\end{theorem}

\section{Holomorphic symmetric embeddings of symmetric domains}
\subsection{Symmetric subdomains}
We continue with the notations of the previous paragraph. Hence $G$ is a
real Lie group of hermitian type (reductive), $\cD$ is the corresponding
domain. We wish to consider reductive subgroups $N\inn G$, also of
hermitian type, defining domains $\cD_N$, such that the inclusion $N\inn G$
induces a holomorphic injection of the domains $i:\cD_N\inn \cD$, and the
$i(\cD_N)$ are totally geodesic. First of all we may assume that $K_N$, a
maximal compact subgroup of $N$, is the intersection $K_N=K\cap N$;
equivalently, letting $o\in \cD$ and $o_N\in \cD_N$ denote the base points,
$i(o_N)=o$. Note that conjugating $N$ by an element of $K$ yields an
isomorphic group $N'$ and subdomain $i':\cD_{N'}\inn \cD$ such that
$i'(o_{N'})=o$, and this defines an equivalence relation on the set of
reductive subgroups $N\inn G$ as described. For the irreducible hermitian
symmetric domains, the equivalence classes of all such $N$ have been
determined by Satake and Ihara (\cite{S1} for the cases of $\cD$ of type
$\bf I_{p,q},\ II_n, \ III_n$; \cite{I} for the other cases).

Before quoting the results we will need, let us briefly remark on the
mathematical formulation of the conditions. For this, let $\cD,\ \cD'$ be
hermitian symmetric domains, $G,\ G'$ the automorphism groups, $\Gg,\ \Gg'$
the Lie algebras, $\Gg=\kk\oplus\pp,\ \Gg'=\kk'\oplus\pp'$ the Cartan
decompositions and $\gt,\ \gt'$ the Cartan involutions on $\Gg$ and $\Gg'$,
respectively. To say that for an injection $i_{\cD}:\cD\hra \cD'$ of
symmetric spaces,
$i_{\cD}(\cD)$ is
totally geodesic in $\cD'$ is to say that
$i_{\cD}$ is induced by an injection $i:\Gg\hra \Gg'$ of the Lie
algebras. If this holds, $i_{\cD}$ is said to be {\it strongly
  equivariant}. Then, $\gt=\gt'_{|i(\Gg)}$, or $\kk=\Gg\cap \kk',\
\pp=\Gg\cap \pp'$. Since both $\cD$ and $\cD'$ are hermitian
symmetric, there is an element $\xi$ in the center of $\kk$ (resp.
$\xi'$ in the
center of $\kk'$), such that $J=\ad(\xi)$ (resp. $J'=\ad(\xi')$) gives the
complex structure. To say that the injection $i_{\cD}:\cD\hra \cD'$ is
holomorphic is the same as saying $i\circ J = J'\circ i$, or equivalently,

\vspace*{.2cm}
$\hbox{(H$_1$)}\hspace*{5.8cm} i\circ \ad(\xi) = \ad(\xi')\circ i.$

\vspace*{.2cm}
\noindent This is the condition utilized by Satake and Ihara in their
classifications. The condition (H$_1$) is clearly implied by

\vspace*{.2cm}(H$_2$)\hspace*{6.5cm}$i(\xi)=\xi',$

\vspace*{.2cm}
\noindent which however, if fulfilled, gives additional information. For
example (\cite{S2}, Proposition 4) if $\cD$ is a tube domain and $i$
satisfies (H$_2$), then $\cD'$ is also a tube domain.
Furthermore, (\cite{S}, Proposition II 8.1), if $i_{\cD}:\cD\lra \cD'$
is a holomorphic map which is strongly equivariant, the corresponding
homomorphism $i$
fulfills (H$_1$), and, moreover, if $\cD$ and $\cD'$ are viewed as bounded
symmetric domains ${\bf D}$, $\bf D'$ via the Harish-Chandra embeddings,
then $i_{\hbf{D}}:{\bf D}\lra {\bf D}'$
is the restriction of a $\fC$-linear map
$i^+:\pp^+\lra (\pp')^+$. If $i_{\fC}:\Gg_{\fC}\lra \Gg_{\fC}'$ denotes the
$\fC$-linear extension of $i$, and $\gs:\Gg_{\fC}\lra \Gg_{\fC},\
\gs':\Gg_{\fC}'\lra \Gg_{\fC}'$ denote the conjugations over $\Gg$ and
$\Gg'$, respectively, then the condition $\gt=\gt'_{|i(\Gg)}$ is equivalent
to the condition $i_{\fC}\circ \gs=\gs'\circ i_{\fC}$. This implies that
$i:(\Gg,\xi)\lra (\Gg',\xi')$ gives rise to a symmetric Lie algebra
homomorphism $(\Gg_{\fC},\gs)\lra (\Gg_{\fC}',\gs')$, and therefore, by
\cite{S}, Proposition I 9.1, to a homomorphism of Jordan triple systems
$i^+:\pp^+\lra (\pp')^+$. It follows (\cite{S}, p.~85) that the following
three categories are equivalent:
\begin{tabular}{cc}$(\scS\scD)$ & \parbox[t]{15cm}{Category whose objects
  are symmetric
  domains $(\cD,o)$ with base point $o$, whose morphisms
  $\gr_{\cD}:(\cD,o)\lra (\cD',o')$ are strongly equivariant holomorphic
  maps $\gr_{\cD}:\cD\lra \cD'$ with $\gr_{\cD}(o)=o'$.} \\
   $(\scH\scL)$ & \parbox[t]{15cm}{Category whose objects are
  semisimple Lie algebras
  $(\Gg,\xi)$ of hermitian type (without compact factors), whose morphisms
  $\gr:(\Gg,\xi)\lra (\Gg',\xi')$ are homomorphisms satisfying (H$_1$).} \\
   $(\scH\scJ)$ & \parbox[t]{15cm}{Category whose objects are
  positive definite hermitian
  Jordan triple systems $\pp^+$, whose morphisms $\gr_+:\pp^+\lra (\pp')^+$
  are $\fC$-linear homomorphisms of Jordan triple systems.}
\end{tabular}

\subsection{Positive-dimensional boundary components}
We now quote some results which we will be using. First, assume we have
fixed $A\inn G$ as above, and let $F_b\inn
\overline{\bf D}$ be a standard boundary component of positive dimension,
i.e., if $\cD$ is irreducible, of rank $t$, then $b<t$; if $\cD=\cD_1\times
\cdots\times \cD_d$, then in the notations of (\ref{e3.2}), ${\bf
  b}=(b_1,\ldots, b_d)$, we have $b_i<t_i$ for {\it at least} one $i$. If
$\cD$ is irreducible, we list in Table \ref{T1} a positive-dimensional
boundary component and a symmetric subdomain $\cD_M\inn \cD$ with the
property that $\cD_N=\cD_F\times \cD'$, where $\cD_F$ is, as a hermitian
symmetric space, isomorphic to the given boundary component. If $\cD$ is
reducible, $\cD=\cD_1\times\cdots\times\cD_d$, and $F_{\bfb}\inn
\overline{\bf D}$ is a standard boundary component, we get a subdomain
$\cD_N=\cD_{N_1}\times \cdots \times \cD_{N_d}$ such that $\cD_{N_i}\inn
\cD_i$ is of the type just mentioned with respect to the boundary component
$F_{b_i}\inn \overline{\bf D}_i$.

\begin{table}\caption{\label{T1} Symmetric subdomains incident with
    positive-dimensional boundary components}
  $$\begin{array}{|c|c|c|c|}\hline \cD & F_b,\ \ (b<t) & \cD_N &
    \hbox{(H$_2$)} \\ \hline \hline \bf I_{\hbox{\scsi\bf p,q}} & \bf
    I_{\hbox{\scsi\bf p-b,q-b}} & \bf I_{\hbox{\scsi\bf p-b,q-b}}\times
    I_{\hbox{\scsi\bf b,b}} & p=q \\ \hline \bf II_{\hbox{\scsi\bf n}} & \bf
    II_{\hbox{\scsi\bf n-2b}} & \bf II_{\hbox{\scsi\bf n-2b}}\times
    II_{\hbox{\scsi\bf 2b}} & yes \\ \hline \bf III_{\hbox{\scsi\bf n}} & \bf
    III_{\hbox{\scsi\bf n-b}} & \bf III_{\hbox{\scsi\bf n-b}}\times
    III_{\hbox{\scsi\bf b}} & yes \\ \hline \bf IV_{\hbox{\scsi\bf n}} & \bf
    IV_{\hbox{\scsi\bf 1}} & \bf IV_{\hbox{\scsi\bf 1}}\times
    IV_{\hbox{\scsi\bf 1}} & yes \\ \hline \bf V & \bf I_{\hbox{\scsi\bf
        5,1}} & \bf I_{\hbox{\scsi\bf 5,1}}\times I_{\hbox{\scsi\bf 1,1}} &
    yes \\ \hline \bf VI & \bf IV_{\hbox{\scsi\bf 10}} & \bf
    IV_{\hbox{\scsi\bf 10}}\times IV_{\hbox{\scsi\bf 1}} & yes \\ & \bf
    IV_{\hbox{\scsi\bf 1}} & \bf IV_{\hbox{\scsi\bf 1}}\times
    IV_{\hbox{\scsi\bf 10}} & yes \\ \hline
\end{array}$$
\end{table}
Next, choose $N\inn G$ with $\cD_N$
 as in Table \ref{T1}, such that $A\inn N$ is a
 maximal $\fR$-split torus in $N$, so that we can speak of standard boundary
 components of $\cD_N$. Then the subdomains $\cD_N$ listed in Table
 \ref{T1} have the following property. For simplicity we will assume from
 now on that $G$ is semisimple.
\begin{proposition}\label{p7.1} Given $G$, simple of hermitian type with
  maximal $\fR$-split torus $A$ and simple $\fR$-roots $\eta_i$ ($1\leq
  i\leq t=\rank_{\fR}G$), let $P_{b}$ and $F_{b}$ denote the
  standard maximal $\fR$-parabolic and standard boundary component
  determined by $\eta_b$ ($b<t$). Let $N\inn G$ be a symmetric subgroup
  with $A\inn N$, defining a subdomain $\cD_N$ as in Table \ref{T1}, such
  that $N=N_{1}\times N_{2}$ and $N_{1}$ is a hermitian Levi factor
  of $P_b$.
  Let $P_0\times P_{{t_2}}$ be the standard maximal parabolic defined
  by the last simple $\fR$-root $\eta_{t_2}$ of the second factor in the
  decomposition
  $N=N_{1}\times N_{2}$. Then if $F:=F_0\times F_{{t_2}}$
  ($\cong \cD_{N_{1}}\times \{pt.\}$, $t_2=\rank_{\fR}N_{2}$) denotes
  the corresponding standard boundary component, the equality
  $i_N(F)=F_{b}$ holds,
  where $i_N:\cD_{N}\lra \cD$ denotes the injection.
\end{proposition}
{\bf Proof:} From construction, $F\cong\cD_{N_{1}}\times \{pt\}
\cong F_{b}$ as a
hermitian symmetric space; to see that they coincide under $i_N$, recall from
(\ref{e2.3}) the root space decomposition of the hermitian Levi component
of $P_{b}$. Since $A\inn N$ is also a maximal $\fR$-split torus of
$N$, in the root system $\Phi(A,N)$ we have the subsystem
$[\eta_{b+1},\ldots, \eta_t]$ giving rise, on the one hand to the hermitian
Levi factor $\ll_b$ in $\pp_{b}$, on the other hand to the Lie algebra
of the first factor $\nn_{1}$ of $N$. From this it follows that
$P_{b}$ stabilizes $i_N(F)$, hence $i_N(F)=F_{b}$. \ende
Before proceeding to the case of zero-dimensional boundary components, we
briefly explain how the subgroups $N\inn G$ (which are not unique, of
course) arise in terms of $\pm$symmetric/hermitian forms, at least for the
classical cases. For this, we note that $G$ can be described as follows (we
describe here certain reductive groups; the simple groups are just the
derived groups):
\begin{equation}\label{e7.1}\begin{minipage}{14.5cm}\begin{itemize}\item[I]
      $\bf I_{\hbox{\scsi\bf p,q}}$:
      $G$ is the unitary group of a hermitian form on
      $\fC^{p+q}$ of signature $(p,q)$ ($p\geq q$).
\item[II] $\bf II_{\hbox{\scsi\bf n}}$:
      $G$ is the unitary group of a skew-hermitian form on $\fH^n$.
\item[III] $\bf III_{\bfn}$: $G$ is the unitary group of a skew-symmetric form
     on $\fR^{2n}$.
\item[IV] $\bf IV_{\bfn}$:
     $G$ is the unitary group of a symmetric bilinear form on
     $\fR^{n+2}$ of signature $(n,2)$.
\end{itemize}\end{minipage}\end{equation}
Each of the $\pm$symmetric/hermitian forms is isotropic, and if
$t=\rank_{\fR}G$, the maximal dimension of a totally isotropic subspace is
$t=q,\ \left[{n\over 2}\right],\ n,\ 2$ in the cases I, II, III, and IV,
respectively. Each maximal real parabolic is the stabilizer of a
totally isotropic
subspace, and using the canonical order on the $\fR$-roots as above,
$P_{b}$ stabilizes a totally isotropic subspace of dimension
$b$. Choosing a maximal torus $T$ (resp. a maximal $\fR$-split torus
$A\inn T$) amounts to choosing a basis of $V$ (resp. choosing a subset of
this basis which spans a
maximal totally isotropic subspace), and the standard parabolic is the
stabilizer of a totally isotropic subspace spanned by some part of this
basis. Now let $H\inn V$ be a totally isotropic subspace with basis
$h_1,\ldots, h_b$. Then there are elements $h_1',\ldots,h_b'$ of $V$ such
that $h(h_i,h_j')=\gd_{ij},\ h(h_i,h_i)=h(h_i',h_i')=0$, and $h_1,\ldots,
h_b,h_1',\ldots, h_b'$ span (over $D$) a vector subspace $W\inn V$ on which $h$
restricts to a non-degenerate form. Let $W^{\perp}$ denote the orthogonal
complement of $W$ in $V$, $W\oplus W^{\perp}=V$. Then
\begin{equation}\label{e8.1} N:=U(W,W^{\perp};h):=\{g\in U(V,h) |
  g(W)\inn W, g(W^{\perp})\inn W^{\perp} \}\cong U(W,h_{|W})\times
  U(W^{\perp},h_{|{W^{\perp}}}).
\end{equation}
$N$ is a reductive subgroup of $G$, and as one easily sees, its symmetric space
is just the domain denoted $\cD_N$ in Table \ref{T1} above.
The relation
``boundary component $\inn $ symmetric subdomain'' translates into
``totally isotropic subspace $\inn$ non-degenerate subspace'', $H\inn W$,
and {\it because} $h_{|W}$ is non-degenerate, any $g\in U(V,h)$ which
stabilizes $W$ automatically stabilizes its orthogonal complement in $V$ as
in (\ref{e8.1}).

\subsection{Zero-dimensional boundary components}
We now would like to consider the zero-dimensional boundary components,
which correspond in the above picture to maximal totally isotropic
subspaces.
The construction above (\ref{e8.1}) doesn't necessarily work in this case,
as $W^{\perp}$ may be $\{0\}$, and $N=G$.
However, in terms of
domains, given {\it any} subdomain $\cD'\inn \cD$, it can be translated so
as to contain a given zero-dimensional boundary component. We therefore
place the following three conditions on such a subdomain:
\begin{itemize}\item[1)] The subdomain $\cD'$ has maximal rank
  ($\rank_{\fR}G'=\rank_{\fR}G$).
\item[2)] The subdomain $\cD'$ is maximal and $G'$ is a maximal subgroup, or
\item[2')] The subdomain $\cD'$ is maximal of tube type and $G'$ is maximal
  with this property.
\item[3')] The subdomain $\cD'$ is maximal irreducible, and $F$ is a
  boundary component of $\cD'$.
\end{itemize}
In Table \ref{T2} we list the subdomains (after \cite{I}) $\cD'$ fulfilling
1), 2) and 3') in the column titled ``$\cD_N$''.
We have listed also those subgroups fulfilling 1),
2') and 3') in the column titled ``maximal tube''.

\begin{table}\caption{\label{T2} Symmetric subdomains incident with
    zero-dimensional boundary components}
$$\begin{array}{|c|c|c|c|} \hline
\cD & \cD_N & \hbox{(H$_2$)} & \hbox{maximal tube} \\ \hline \hline
{\bf I_{\hbox{\scsi\bf p,q}}},\ p>q & \bf I_{\hbox{\scsi\bf p-1,q}} & no &
\bf I_{\hbox{\scsi\bf q,q}}
\\ \hline
\bf I_{\hbox{\scsi\bf q,q}} & - & - & -  \\ \hline
\bf II_{\hbox{\scsi\bf n}},\ n\hbox{ even} & - & - & - \\ \hline
\bf II_{\hbox{\scsi\bf n}},\ n\hbox{ odd} & \bf II_{\hbox{\scsi\bf n-1}}
  & yes & \bf II_{\hbox{\scsi\bf n-1}} \\ \hline
\bf III_{\hbox{\scsi\bf n}} & - & - & - \\ \hline
\bf IV_{\hbox{\scsi\bf n}} & \bf IV_{\hbox{\scsi\bf n-1}} & yes & \bf
IV_{\hbox{\scsi\bf n-1}} \\ \hline
\bf V & \bf I_{\hbox{\scsi\bf 2,4}},\ II_{\hbox{\scsi\bf 5}},\
IV_{\hbox{\scsi\bf 8}} & yes,\ no,\ no & \bf I_{\hbf{2,2}}, II_{\hbf{4}},
    IV_{\hbox{\scsi\bf 8}} \\ \hline
\bf VI & \bf I_{\hbox{\scsi\bf 3,3}},\ II_{\hbox{\scsi\bf 6}} & yes & \bf
I_{\hbox{\scsi\bf 3,3}},\ II_{\hbox{\scsi\bf 6}} \\ \hline
\end{array}$$
{\small In the column ``$\cD_N$'' the subgroups fulfilling 1), 2) and 3')
  are listed, and in the column ``maximal tube'' the subgroups fulfilling
  1), 2') and 3') (i.e., not necessarily 2)) are listed.}
\end{table}

In Table \ref{T2}, if there is no entry in the column
``$\cD_N$'', no such subgroups exist. In these cases it is natural to take the
polydisc $\cD_{N_{\Psi}}$ defined by the maximal set of
strongly orthogonal roots $\Psi=\{\pm\{\mu_1\},\ldots,\pm\{\mu_t\}\}$
as the subdomain
$\cD_N$, as there is no irreducible subdomain, and other products already
occur in Table \ref{T1}. Hence for these cases we require the conditions
2'') and 3'') of the introduction. To sum up these facts we make the following
definition.
\begin{definition}\label{d9.1} Let $G$ be a simple real Lie group of
  hermitian type, $A$ a fixed maximal
  $\fR$-split torus defined as above by a maximal set of strongly
  orthogonal roots, $\eta_i,\ i=1,\ldots, t$ the simple $\fR$-roots,
  $F_{b}$ a standard boundary component and
  $P_{b}$ the corresponding standard maximal $\fR$-parabolic. A
  reductive subgroup $N\inn G$ (respectively the subdomain $\cD_{N}\inn
  \cD$) will be called {\it
    incident} to $P_{b}$ (respectively to $F_{b}$), if $\cD_N$ is
  isomorphic to the corresponding domain of Table \ref{T1} ($b<t$) or Table
  \ref{T2} ($b=t$), and if $N$ fulfills:
  \begin{itemize}\item $b<t$, then $N$ satisfies 1), 2), 3).
   \item $b=t,\ \cD\not\in (\cE\cD)$, then $N$ satisfies 1), 2'), 3').
   \item $b=t,\ \cD\in (\cE\cD)$, then $N$ satisfies 1), 2''), 3'').
  \end{itemize}
  For reducible $\cD=\cD_1\times \cdots \times \cD_d$, we have the product
  subgroups $N_{b_1,1}\times \cdots \times N_{b_d,d}$, where
  $\cD_{N_{b_i,i}}$ is
  incident to the standard boundary component $F_{{b_i}}$ of $\cD_i$
  (and $N_{0,i}=G_i$).
\end{definition}
Next we briefly discuss uniqueness. We consider first the case of
positive-dimensional boundary components. Let $P_{b}$, $1\leq b< t$
be a standard parabolic and let $L_b$ be the ``standard'' hermitian Levi
factor, i.e., such that $Lie(L_b)=\ll_b$; then
\begin{equation}\label{e10.1} N_b:= L_b\times \cZ_G(L_b)
\end{equation}
is a subgroup having the properties of Proposition \ref{p7.1}, unique since
$L_b$ is unique. We shall refer to this unique subgroup as the {\it
  standard} incident subgroup. The different Levi factors
$\cL$ in Levi decompositions
$P_b=\cL\rtimes \cR_u(P_b)$ are conjugate by elements of $\cR_u(P_b)$, as is
well known. This implies for the hermitian factors
$L=\cL^{herm}\inn \cL$ (which
are uniquely determined by $\cL$) by Theorem \ref{t4.1} the following.
\begin{lemma}\label{l10A} Two hermitian Levi factors $L,\ L'\inn P_b$ are
  conjugate by an element of $V_b\inn P_b$.
\end{lemma}
It follows, since $g(L_b\times \cZ_G(L_b))g^{-1}=gL_bg^{-1}\times
\cZ_G(gL_bg^{-1})$, that two subgroups $N,\ N'$, both incident with
$P_b$, are conjugate by an element of $V_b$:
\[\hbox{$N,\ N'$ incident to $P_b$ $\iff$ $N,\ N'$ conjugate (in $G$)
  by $g\in V_b$.}\]
\begin{proposition} If $(N,P_b)$ are incident, there is $g\in V_b$ such
  that $N$ is conjugate by $g$ to the standard $N_b$ of (\ref{e10.1}).
\end{proposition}
{\bf Proof:} Since $N$ is incident, $N\cong N_1\times N_2$, where $N_1$ is
a hermitian Levi factor of $N$. By Lemma \ref{l10A}, $N_1$ is conjugate by
$g\in V_b$ to $L_b$, the hermitian Levi factor with Lie algebra $\ll_b$
in the notations of the last section. Hence $gNg^{-1}=g(N_1\times
N_2)g^{-1} = gN_1g^{-1}\times gN_2g^{-1}=L_b\times N_{2,b}=N_b$, with
$N_{2,b}=\cZ(L_b)$ (this follows from the maximality of $N_b$).
Consequently, $N$ is conjugate by $g\in V_b$ to $N_b$. \ende

The situation for zero-dimensional boundary components is more complicated,
so we just observe the following. Suppose $\cD\not\in (\cE\cD)$, and that
$\cD_N\inn \cD$ is incident to $F_t$, $F_t$=point. For any $g\in
N(F_t)=P_t,\ g\cD_N=\cD'\inn \cD$ is another subdomain, again incident to
$F_t$. If $g\in P_t\cap N_t$, then $g\cD'=\cD_N$. In this sense, letting
$Q_t=P_t\cap N_t$, the coset space $P_t/Q_t$ is a parameter space of
subdomains incident with $F_t$.

Above we have defined the notion of symmetric subgroups incident with a
standard parabolic. Any maximal $\fR$-parabolic is conjugate to one and
only one standard maximal parabolic, $P=gP_bg^{-1}$ for some $b$. Let $N_b$
be any symmetric subgroup incident with $P_b$. Then just as one has the
pair $(P_b,N_b)$ one has the pair $(P,N)$,
\begin{equation}\label{e10.3} P=gP_bg^{-1},\quad N=g N_b g^{-1}.
\end{equation}
\begin{definition} \label{d10.1} A pair $(P,N)$ consisting of a maximal
  $\fR$-parabolic $P$ and a symmetric subgroup $N$ is called {\it
    incident}, if the groups are conjugate by a common element $g$ as in
  (\ref{e10.3}) to a pair
  $(P_b,N_b)$ which is incident as in Definition \ref{d9.1}.
\end{definition}

\section{Rational parabolic and rational symmetric subgroups}
\subsection{Notations}
We now fix some notations to be in effect for the rest of the paper. We
will be dealing with algebraic groups defined over $\fQ$, which give rise
to hermitian symmetric spaces, groups of {\it hermitian type}, as we will
say. As we are interested in the automorphism groups of domains, we may,
without restricting generality, assume the group is {\it centerless}, and
{\it simple} over $\fQ$. Henceforth $G$ will denote such an algebraic
group. To avoid complications, we exclude in this paper the following case:

\vspace*{.2cm}{\bf Exclude:}\hspace*{2cm} All non-compact real factors of
$G(\fR)$ are of type $SL_2(\fR)$.

\vspace*{.2cm}
\noindent Finally, we shall
only consider {\it isotropic} groups. This implies the hermitian symmetric
space $\cD$ has no compact factors. By our assumptions, then, we have
\begin{itemize}\item[(i)] $G=Res_{k|\fQ}G'$, $k$ a totally real number
  field, $G'$ absolutely simple over $k$.
\item[(ii)] $\cD=\cD_1\times \cdots \times \cD_d$, each $\cD_i$ a non-compact
  irreducible hermitian symmetric space, $d=[k:\fQ]$.
\end{itemize}
We now introduce a few notations concerning the root systems involved. Let
$\gS_{\infty}$ denote the set of embeddings $\gs:k\lra \fR$; this set is in
bijective correspondence with the set of infinite places of $k$. We
denote the latter by $\nu$, and if necessary we denote the
corresponding embedding by
$\gs_{\nu}$. For each $\gs\in \gS_{\infty}$, the group $^{\gs}G'$ is the
algebraic group defined over $\gs(k)$ by taking the
set of elements $g^{\gs},\ g\in G'$. For each infinite
prime $\nu$ we have
$G_{k_{\nu}}\cong (^{\gs_{\nu}}G')_{\fR}$, and the decomposition of $\cD$
above can be written
$$\cD=\prod_{\gs\in \gS_{\infty}}\cD_{\gs},\quad
\cD_{\gs}:=(^{\gs}G')_{\fR}/K_{(\gs)}=(^{\gs}G')_{\fR}^0/K_{(\gs)}^0.$$
Since $G'$ is isotropic, there is a positive-dimensional $k$-split torus
$S'\inn G'$, which we fix. Then  ${^{\gs}S}'$ is a maximal $\gs(k)$-split
torus of $^{\gs}G'$ and there is a canonical isomorphism $S'\ra {^{\gs}S}'$
inducing an isomorphism $\Phi_k=\Phi(S',G')\lra
\Phi_{\gs(k)}({^{\gs}S}',{^{\gs}G}')=:\Phi_{k,\gs}$.
The torus $Res_{k|\fQ}S'$ is
defined over $\fQ$ and contains $S$ as maximal $\fQ$-split torus; in fact
$S\cong S'$, diagonally embedded in $Res_{k|\fQ}S'$. This yields an
isomorphism $\Phi(S,G)\cong \Phi_k$, and the root systems
$\Phi_{\fQ}=\Phi(S,G)$, $\Phi_k$ and $\Phi_{k,\gs}$ (for all $\gs\in
\gS_{\infty}$) are identified
by means of the isomorphisms.

In each group $^{\gs}G'$ one chooses a maximal $\fR$-split torus
$A_{\gs}\nni {^{\gs}S}'$, contained in a maximal torus defined over
$\gs(k)$. Fixing an order on $X(S')$ induces one also on $X({^{\gs}S}')$
and $X(S)$. Then, for each $\gs$, one chooses an order on $X(A_{\gs})$
which is compatible with that on $X({^{\gs}S}')$, and $r:X(A_{\gs})\lra
X({^{\gs}S}')\cong X(S)$ denotes the restriction homomorphism. The canonical
numbering on $\gD_{\fR,\gs}$ of simple $\fR$-roots of $G$ with respect to
$A_{\gs}$ is compatible by restriction with the canonical numbering of
$\gD_{\fQ}$ (\cite{BB}, 2.8). Recall also that each $k$-root in $\Phi_k$ is
the restriction of at most one simple $\fR$-root of $G'(\fR)$ (which is a
simple Lie group). Let $\gD_k=\{\gb_1,\ldots,\gb_s\}$; for $1\leq i\leq s$ set
$c(i,\gs)$:= index of the simple $\fR$-root of $^{\gs}G'$ restricting on
$\gb_i$. Then $i<j$ implies $c(i,\gs)<c(j,\gs)$ for all $\gs\in \gS$.

Each simple $k$-root defines a unique standard boundary component: for
$b\in \{1,\ldots,s\}$,
\begin{equation}\label{e9.1} F_{\bfb}:=\prod_{\gs\in
    \gS_{\infty}}F_{c(b,\gs)},
\end{equation}
which is the product of standard (with respect to $A_{\gs}$ and
$\gD_{\fR,\gs}$) boundary components $F_{c(b,\gs)}$ of $\cD_{\gs}$. It
follows that $\overline{F}_{\hbox{\scsi\bf j}}\inn
\overline{F}_{\hbox{\scsi\bf i}}$ for $1\leq i\leq j\leq
s$. Furthermore, setting $o_{\bfb}:=\prod o_{c(b,\gs)}$, then (\cite{BB},
p.~472)
\begin{equation}\label{e9.2} F_{\bfb}=o_{\bfb}\cdot L_{\bfb},
\end{equation}
where $L_{\bfb}$ denotes the hermitian Levi component (\ref{e3.3})
of the parabolic
$P_{\bfb}(\fR)=N(F_{\bfb})$. As these are the only boundary components of
interest to us, we will henceforth refer to any conjugates of the
$F_{\bfb}$ of (\ref{e9.1}) by elements of $G$ as {\it rational boundary
  components} (these should more precisely be called rational with respect
to $G$), and to the conjugates of the parabolics $P_{\bfb}$ as the {\it
  rational parabolics}.

\subsection{Rational parabolics}
Let $G'$ be as above, $\gD_k=\{\gb_1\ldots,\gb_s\}$ the set of simple
$k$-roots (having fixed a maximal $k$-split torus $S'$ and an order on
$X(S')$). For $b\in \{1,\ldots,s\}$ we have the standard maximal
$k$-parabolic $P_b'$ of $G'$, whose group of $\fR$-points is the
normalizer of the standard rational boundary component $F_{c(b)}$ of the domain
$\cD'=G_{\fR}'/K'$, where $c(b)$ denotes the index of the simple $\fR$-root
restricting to $\gb_b$; since $G'$ is absolutely simple, $G_{\fR}'$ is simple
and $\cD'$ is irreducible. Hence Theorem
\ref{t4.1} applies to $P_b'(\fR)$. Of the factors given there,
the following are defined over $k$:
the product $M_b'L_b'$ as well as $L_b'$ (but $M_b'$ is {\it not} defined
over $k$, so the $k$-subgroups are (instead of $L_b'$ and $M_b'$) $L_b'$
and ${G_b'}^{(1)}:=M_b'L_b'$), $\cR_b',\cZ_b'$ and $V_b'$. As is well
known, any maximal $k$-parabolic of $G'$ is conjugate to one and only one
of the $P_b'$, and two parabolics are conjugate $\iff$
they are conjugate over $k$. There is a 1-1 correspondence between the set
of $k$-parabolics of $G'$ and the set of $\fQ$-parabolics of $G$, given by
$P'\mapsto Res_{k|\fQ}P'=:P$. The standard maximal $\fQ$-parabolic
$P_{\hbox{\scsi\bf b}}$ of $G$ gives a $\fQ$-structure on the real parabolic
$P_{\hbox{\scsi\bf b}}(\fR)$, which is the normalizer in $\cD$ of the
standard boundary component $F_{\bfb}$ as in (\ref{e9.1}) (see also
(\ref{e3.2}) and (\ref{e3.3})), where ${\bf
  b}=(c(b,\gs_1),\ldots,c(b,\gs_d))$. In the decomposition of Theorem
\ref{t4.1}, the factors $G_{\hbox{\scsi\bf b}}^{(1)}=M_{\hbox{\scsi\bf
    b}}L_{\hbox{\scsi\bf b}},\ \cR_{\hbox{\scsi\bf b}}, \cZ_{\hbox{\scsi\bf
    b}}$ and $V_{\hbox{\scsi\bf b}}$ are all defined over $\fQ$. In
particular, for the factor $G_{\hbox{\scsi\bf b}}^{(1)}$, which we will
call the $\fQ$-hermitian Levi factor (and similarly, we will call
${G_{b}'}^{(1)}$
the $k$-hermitian Levi factor of $P_{b}'$), we have
\begin{equation}\label{e12.1} G_{\hbox{\scsi\bf b}}^{(1)}(\fQ)\cong
  \prod_{\gs}({{^{\gs}G}'}_b^{(1)})_{\gs(k)},\quad \cZ_G(G_{\hbox{\scsi\bf
      b}}^{(1)})(\fQ)\cong \prod_{\gs}(\cZ_{(^{\gs}G'_{\gs(k)})}(
      ({{^{\gs}G}'_b}^{(1)})_{\gs(k)}).
\end{equation}
Furthermore, the hermitian Levi factor $L_{\bfb}$ is defined over $\fQ$,
and
\[L_{\bfb}(\fQ)=\prod_{\gs}({^{\gs}L}_b')_{\gs(k)}. \]

We now make a few remarks about the factors of $G(\fR)$ and of
$L_{\bfb}(\fR)$. Since the map $G'\lra {^{\gs}G}'$ is an isomorphism of a
$k$-group onto a $\gs(k)$-group, the algebraic groups (over $\fC$) are
isomorphic, hence the various ${^{\gs}G}'_{\fR}$ are all $\fR$-forms of
some fixed algebraic group. Similarly, the factors of $L_{\bfb}(\fR)$ are
all $\fR$-forms of a single $\fC$-group. However, they need not be
isomorphic, unless the given $\fC$-group has a unique $\fR$-form of
hermitian type (like $Sp(2n,\fC)$). Next we note the following.
\begin{lemma}\label{L12a} $L_{\bfb}$ is anisotropic $\iff$ $b=s$.
\end{lemma}
{\bf Proof:} The group $L_{\bfb}$ is anisotropic precisely when the
boundary component $F_{\bfb}$ defined by it contains no other boundary
components $F_{\hbf{c}}^*\inn F_{\bfb}^*$, which means $b\geq c$ for all
$c$, or $b=s$. \ende
In this case the group $L_{\bfb}$ does not fulfill the assumptions we have
placed on $G$, and our results up to this point are not directly applicable
to $L_{\bfb}$. Let us see how the phenomenon of compact factors of
$L_{\bfb}(\fR)$ manifests itself in
$F_{\bfb}=\prod_{\gs}F_{c(b,\gs)}$. Suppose some factor of $L_{\bfb}(\fR)$
is compact, say $L_{1,b}$. Then the symmetric space $\cD_{b,\gs_1}$
of $L_{1,b}$ is compact, so it is not true that $\cD_{b,\gs_1}\cong
F_{c(b,\gs_1)}$, hence it is also not true that $\cD_{\bfb}\cong F_{\bfb}$,
where $\cD_{\bfb}=\prod_{\gs}\cD_{b,\gs}$ is the symmetric space of
$L_{\bfb}(\fR)$. However, letting $\cD_{\bfb}'$ be the product of all
compact factors, $\cD_{\bfb}/\cD_{\bfb}'$ is a symmetric space which
is isomorphic to $F_{\bfb}$. What happens is that in the product
$F_{\bfb}=\prod F_{c(b,\gs)}$, all factors $F_{c(b,\gs)}$ are {\it points}
for which $\cD_{b,\gs}$ is {\it compact}. Hence whether this occurs depends
on whether any factors $\cD_{\gs}$ have zero-dimensional (rational)
boundary components or not.

\subsection{Incidence}
We keep the notations used above; $G$ is a simple $\fQ$-group of hermitian
type. Our main definition gives a $\fQ$-form of Definition \ref{d10.1}, and
is the following.
\begin{definition}\label{d12.1} Let $P\inn G$ be a maximal $\fQ$-parabolic,
  $N\inn G$ a reductive $\fQ$-subgroup. Then we shall say that $(P,N)$ are
  {\it incident} (over $\fQ$), if $(P(\fR),N(\fR))$ are incident in the
  sense of Definition \ref{d10.1}.
\end{definition}
Note that in particular $N$ must itself be of hermitian type, and such that
the Cartan involution of $G(\fR)$ restricts to the Cartan involution of
$N(\fR)$. Furthermore, $N$ must be a $\fQ$-form of a product of groups,
defining domains each of which is as in either Table \ref{T1}
or Table \ref{T2}.

The main result of this paper is the following existence result.

\begin{theorem}\label{t12.1} Let $G$ be $\fQ$-simple of hermitian type
  subject to the restrictions above ($G$ is isotropic and $G(\fR)$ is not a
  product of $SL_2(\fR)$'s),
  $P\inn G$ a $\fQ$-parabolic. Then there exists a reductive $\fQ$-subgroup
  $N\inn G$ such that $(P,N)$ are incident over $\fQ$, with the exception
  of the indices $C^{(2)}_{2n,n}$ for the zero-dimensional boundary
  components.
\end{theorem}
We will give the proof in the following sections, where we consider separately
different cases (of the $\fQ$-rank, the dimension of a maximal $\fQ$-split
torus). But before we start, we note here that by definition, if the
theorem holds for {\it standard} parabolics, then it holds for all
parabolics, so it will suffice to consider only standard parabolics. The
case that $G'$ has index $C^{(2)}_{2,1}$ was considered in \cite{hyp}; in
that case there is a unique standard parabolic $P_1$, with zero-dimensional
boundary component; the associated $N_1$ described in \cite{hyp}
has domain $\cD_{N_1}$ which
is not a two-disc, but only a one-dimensional disc.

\section{Split over $\fR$ case}
In this paragraph we consider the easiest case. This could loosely be
described as an $\fR$-Chevally form.
\begin{definition}\label{d13.1} Let $G'$ be as in the last paragraph,
  absolutely simple over $k$, and let $\Phi_k$ be a root system
  (irreducible) for $G'$ with respect to a maximal $k$-split torus $S'\inn
  G'$. Let $\Phi_{\fR}$ be the root system of $G'(\fR)$ with respect to a
  maximal $\fR$-split torus $A'$ of the real (simple) group $G'(\fR)$. We
  call $G'$ {\it split over $\fR$}, if $\Phi_k\cong \Phi_{\fR}$ as root
  systems, and if the indices of $G'$ and $G'(\fR)$ coincide.
\end{definition}
Note that the indices are independent of the split tori used to form the
root system, so there is no need to assume $S'\inn A'$ in the above
definition (the notion of isomorphism of indices is obvious). However, one
can always find split tori $S', A'$ such that $S'\inn A'$. From
$\Phi_k\cong \Phi_{\fR}$ it follows then that $S'=A'$, as both tori have
the same dimension.
\begin{lemma}\label{l13.1} Let $G$ be simple over $\fQ$ (=$Res_{k|\fQ}G'$),
  $\cD=\prod_{\gs\in \gS_{\infty}}\cD_{\gs}$ the domain defined by the real
  Lie group $G(\fR)\cong \prod_{\gs\in
    \gS_{\infty}}G_{\gs}=:\prod_{\gs}{^{\gs}G'_{\fR}}$. If $G'$ is split
  over $\fR$, then $G_{\gs}=G_{\tau}$ for all $\gs,\tau\in \gS_{\infty}$.
\end{lemma}
{\bf Proof:} For each $\gs$ we have $A_{\gs}\nni {^{\gs}S}'$, so by
assumption $A_{\gs}\cong {^{\gs}S}'$, and for each $\gs$ the map
$\phi:\Phi_k\lra \Phi_{\gs(k)}(^{\gs}G')$ is an isomorphism, and since
$\Phi_k\cong \Phi_{\fR}$,
$$\Phi_{\gs(k)}(^{\gs}G')\cong \Phi_{\fR}(^{\gs}G_{\fR}').$$ It follows
that $\Phi_{\fR}\cong \Phi_k\stackrel{\phi}{\cong}
\Phi_{\gs(k)}(^{\gs}G')\cong \Phi_{\fR}(^{\gs}G'_{\fR})\cong \Phi_{\fR}$.
Similarly, since the index of $G'$ is isomorphic to the index of $G'(\fR)$
(which determines the isomorphy class of $G'(\fR)$), the index of
${^{\gs}G}'$ is isomorphic to that of ${^{\gs}G}'(\fR)$. But the index of
$G'(\fR)$ is the same as ${^{\gs}G}'(\fR)$, as an easy case by case check
verifies. For example, for type (I), all factors have the same $\fR$-rank
$q$, hence are all isomorphic to $SU(p,q)$. See Examples \ref{examples}
below for the other cases. Hence ${^{\gs}G}'(\fR)\cong {^{\tau}G}'(\fR)$
for all $\gs, \tau$, as claimed. \ende

{}From this it follows in particular that the (standard) boundary components
are determined by $c(b,\gs)=b,\ \forall_{\gs},\ {\bf b}=(b,\ldots,b),\
1\leq b\leq t=\rank_{\fQ}G=\rank_{k}G'=\rank_{\fR}G'(\fR)$. Hence they are
of the form
\begin{equation}\label{e13.1} F_{\bfb}=\prod_{\gs\in
    \gS_{\infty}}F_{b,\gs},
\end{equation}
and $F_{b,\gs}$ is the standard rational boundary component of $\cD_{\gs}$.

\begin{examples}\label{examples}
  We now give examples of split over $\fR$ groups in each of the cases, and
  any such will be of one of the listed types. Let $k$ be a totally real
  number field.
\begin{itemize}\item[I.] Let $K|k$ be imaginary quadratic, $V$ a
  $K$-vector space of dimension $n=p+q$, and $h$ a hermitian form on $V$
  defined over $K$. Then the unitary group $U(V,h)$ is split over $\fR$
  $\iff$ the hermitian form $h$ has Witt index $q$ and for all infinite
  primes, $h_{\nu}$ has signature $(p,q)$.
\item[II.] Let $D|k$ be a totally definite quaternion algebra over $k$
  (with the canonical involution), $V$ an $n$-dimensional right vector
  space over $D$, $h$ a skew-hermitian form on $V$ defined over $k$. Then
  the unitary group $U(V,h)$ is split over $\fR$ $\iff$ the skew-hermitian
  form has Witt index $[{n\over2}]$ ($n>4$).
\item[III.] Take $G=Sp(2n,k)$.
\item[IV.] Let $V$ be a $(n+2)$-dimensional $k$-vector space, $h$ a
  symmetric bilinear form defined over $k$ of Witt index 2. Then if
  $U(V,h)$ is of hermitian type, it is split over $\fR$.
\item[V.] The Lie algebra in this case is of the form
  $\scL(\CC_k,(J_1^b)_k)$, the Tits algebra, where $\CC_k$ is an
  anisotropic octonion algebra and $(J_1^b)_k$ is the Jordan algebra
  $\BB^+$ for an associative algebra $\BB$ whose traceless elements with
  the Lie product form a Lie algebra of type $\Ss\uu(2,1)$; since $G'$ is
  split over $\fR$, the algebra $\BB^-$ is the Lie algebra of a unitary
  group of a $K$-hermitian form ($K|k$ imaginary quadratic as in (I)) of
  Witt index 1.
\item[VI.] The Lie algebra is isomorphic to $\scL(\AA_k,\JJ_k)$, the Tits
  algebra, where $\AA_k$ is a totally indefinite quaternion algebra over
  $k$ and $\JJ_k$ is a $k$-form of the exceptional Jordan algebra denoted
  $J^b$ by Tits.
\end{itemize}
\end{examples}

\begin{lemma}\label{l14.1} In the notations above, let $N'(\fR)\inn
  G'(\fR)$ be a subgroup such that the Lie algebra $\nn'\inn \Gg'$ is a
  {\it regular} subalgebra, i.e., defined by a closed symmetric set of
  roots $\Psi$ of the (absolute) root system $\Phi$ of $G'$. Then $N'$ is
  defined over $k$, $N'\inn G'$.
\end{lemma}
{\bf Proof:} From the isomorphism of the indices of $G'$ and $G'(\fR)$, it
follows that any subalgebra $\Gg'\inn \Gg$, such that for some subset
$\Psi\inn\Phi$, the subalgebra $\Gg'$ is given by $\Gg'=\tt+\sum_{\eta\in
  \Psi}\Gg^{\eta}$, is defined over $\fR$ $\iff$ it is defined over $k$.
The regular subalgebra $\nn'$ is of this type, and it follows that $N'$ is
defined over $k$. \ende
\begin{corollary}\label{c14.1} Let $N'\inn G'$ be as in Lemma \ref{l14.1},
  $N=Res_{k|\fQ}N'\inn G$. Then $N$ is defined over $\fQ$.
\end{corollary}
To apply Lemma \ref{l14.1} to ($k$-forms of) subgroups whose domains are
listed in Tables \ref{T1} and \ref{T2}, we need to know which of the
subgroups are defined by regular subalgebras. Ihara in \cite{I} considered
this question, and the result is: all isomorphism classes of groups in
Table \ref{T1} and all isomorphism classes of groups in Table \ref{T2},
with the exception of $SO(n-1,2)\inn SO(n,2)$ for $n$ even, have
representatives which are defined by (maximal) regular subalgebras.
\begin{corollary}\label{c14.2} Let $G'$ be split over $\fR$. Then Theorem
  \ref{t12.1} holds for $G=Res_{k|\fQ}G'$.
\end{corollary}
{\bf Proof:} By Lemma \ref{l13.1}, $G(\fR)/K=\cD=\prod\cD_{\gs}$, and all
$\cD_{\gs}$ are isomorphic to $\cD'=G'(\fR)/K'$; the rational boundary
components are as in (\ref{e13.1}), products of copies of $F_b'$, the
standard boundary component of $\cD'$, and each $\fQ$-parabolic of $G$ is
conjugate to one of $P_b=Res_{k|\fQ}P_b'$, where $P_b'(\fR)=N(F_b')$ is the
standard maximal real parabolic of $G'(\fR)$. Now locate $F_b'$ in Table
\ref{T1} or \ref{T2} as the case may be; the corresponding group $N_b'$ is
isomorphic to one defined by a regular subalgebra of $\Gg_{\fR}'$ with the
one exception mentioned above. Then by Lemma \ref{l14.1}, $N_b'$ is defined
over $k$, hence (Corollary) $N_b=Res_{k|\fQ}N_b'$ is defined over $\fQ$,
and is incident with $P_b$. This takes care of all cases except the
exception just mentioned, $\bf IV_{\hbox{\scsi\bf n-1}}\inn
IV_{\hbox{\scsi\bf n}}$, $n>3$ even. So let $V$ be a $k$-vector space of
dimension $n+2$, $h$ a symmetric bilinear form on $V$. By assumption, $G'$
is split over $\fR$, so the Witt index of $h$ is 2. Let $H\inn V$ be a
maximal totally isotropic subspace (two-dimensional) defined over $k$, and
$h_1,h_2$ a $k$-basis. Then there are $k$-vectors $h_i'$ such that
$H_1:=<h_1,h_1'>$ and $H_2:=<h_2,h_2'>$ are hyperbolic planes; let
$W=H_1\oplus H_2$ denote their direct sum. From $n>3$, $W$ has codimension
$\geq1$ in $V$. Let $v\in W^{\perp}$ be a $k$-vector, and set:
$$U:=v^{\perp}=\{w\in V|h(v,w)=0\}.$$ Then $W\inn U$, the dimension of $U$
is $n+2-1=n+1$, and $h_{|U}$ still has Witt index 2. Hence
$$N':=\{g\in U(V,h)|g(U)\inn U\}$$ is a $k$-subgroup, and $N'(\fR)^0\isom
SO(n-1,2)$. This is a group which is incident to a parabolic whose group of
real points is the stabilizer of the zero-dimensional boundary component
$F_{2}'$ of the domain $\cD'$ of type $\bf IV_{\hbox{\scsi\bf n}}$. \ende
This completes the discussion of the split over $\fR$ case. We just mention
that, at least in the classical cases, we could have argued case for case
with $\pm$symmetric/hermitian forms as in the proof of the exception above.
Using the root systems simplified the discussion, and, in particular, gives
the desired results for the exceptional groups without knowing their
explicit construction.

\section{Rank $\geq 2$}
In this paragraph we assume $G$ in {\it not} split over $\fR$, but that
$\rank_kG'=\rank_{\fQ}G\geq 2$. Under these circumstances, it is known
precisely which $k$-indices are possible for $G'$ of hermitian type.
\begin{proposition}\label{p16.1} Assume $\rank_{\fQ}G\geq 2$ and that $G'$
  is not split over $\fR$. Then the $k$-index of $G'$ is one of the
  following:
\begin{itemize}\item[(I)] ${^2A}^{(d)}_{n,s};\ s\geq2,\ d|n+1,\ 2sd\leq
  n+1$; if $d=1$, then $2s<n+1$.
\item[(II)] ${^1D}^{(2)}_{n,s},\ s\geq2,\ s< \ell\ (n=2\ell);\quad
  {^2D}^{(2)}_{n,s},\ s\geq 2,\ s<\ell\ (n=2\ell+1)$.
\item[(III)] $C^{(2)}_{n,s},\ s\geq 2,\ s< [{n\over 2}]$.
\item[(IV)] none
\item[(V)] none
\item[(VI)] $E^{31}_{7,2}$.
\end{itemize}
\end{proposition}
{\bf Proof:} All statements are self-evident from the description of the
indices in \cite{T}; in the case (V) there are three possible indices, only
one of which has rank $\geq 2$; this is the split over $\fR$ index.
Similarly, in the case (IV), rank $\geq2$ implies split over $\fR$. For
type (III), the indices $C^{(1)}$ are also split over $\fR$. \ende There is
only one exceptional index to consider, so we start by dealing with this
case. The index we must discuss is
$$\begin{minipage}{16.5cm} \unitlength1.5cm

  
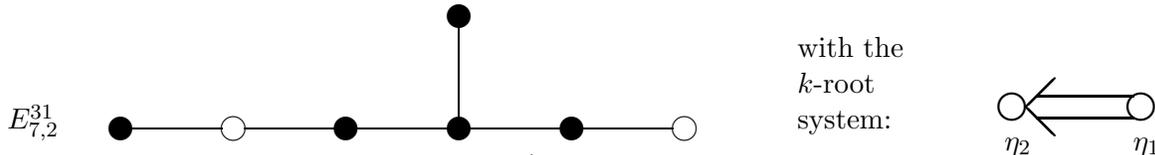
\hspace*{2cm}
\begin{picture}(9,1)
  \put(0,0){\circle*{.2}} \put(1,0){\circle{.2}} \put(2,0){\circle*{.2}}
  \put(3,0){\circle*{.2}} \put(4,0){\circle*{.2}} \put(5,0){\circle{.2}}

  \put(-1,0){$E^{31}_{7,2}$}

  \put(.1,0){\line(1,0){.8}} \put(1.1,0){\line(1,0){.8}}
  \put(2.1,0){\line(1,0){.8}} \put(3.1,0){\line(1,0){.8}}
  \put(4.1,0){\line(1,0){.8}} \put(3,.1){\line(0,1){.8}}

  \put(3,1){\circle*{.2}}

  \put(6,0){\parbox[b]{2cm}{with the \\ $k$-root \\ system:}}

  \put(7,0){\mbox{
\setlength{\unitlength}{0.0037500in}%
\begin{picture}(216,60)(12,710)
  \thicklines \put(320,740){\circle{36}} \put(310,680){$\eta_1$}
  \put(160,740){\line( 1, 1){ 40}} \put(140,740){\circle{36}}
  \put(130,680){$\eta_2$} \put(160,740){\line( 1,-1){ 40}}
  \put(175,725){\line( 1, 0){135}}
%\put(175,755){\makebox(0.4444,0.6667){\SetFigFont{10}{12}{rm}.}}
  \put(175,755){\line( 1, 0){135}}
\end{picture}}}

\end{picture}
\end{minipage}$$
There are two simple $k$-roots, $\eta_1$ and $\eta_2$; let $P_{b}'$ be the
corresponding standard maximal $k$-parabolics, $F_{b}'$ the corresponding
standard boundary components of the irreducible domain $\cD'$. Then
$F_{2}'$ is the one-dimensional boundary component, $F_{1}'$ is the
ten-dimensional one. The $k$-root system is of type ${\bf BC}_2$ (since the
highest simple $\fR$-root is anisotropic, see \cite{BB}, 2.9). Consider the
decomposition of Theorem \ref{t4.1} for $P_{b}'(\fR)$; in both cases $L'_b$
is non-trivial, and, as mentioned above, $M_b'\cdot L_b'$ is defined over
$k$. Here we have $b=1$ or 2. But for $E_7$, the compact factor $M_b'$ is
in fact {\it absent}\footnote{see \cite{S}, p.~117}, and as $L_b'$ is
defined over $k$, we can set
$$N_b'=L_b'\times \cZ_{G'}(L_b').$$ This is a $k$-subgroup which is a
$k$-form of the corresponding $\fR$-subgroup whose domain is listed in
Table \ref{T1}.  Now consider $G=Res_{k|\fQ}G'$. It also has two standard
maximal parabolics $P_{\hbox{\scsi\bf 1}}$ and $P_{\hbox{\scsi\bf 2}}$, and
in each we have a non-trivial hermitian Levi factor\footnote{we note a
  change of notation here in that in (\ref{e3.3}), $L_{\bfb}$ denotes a
  real Lie group} $L_{\bfb}:= Res_{k|\fQ}L_b'$, such that
$$L_{\bfb}(\fR)=\prod_{\gs\in \gS_{\infty}}{^{\gs}(L_b')}_{\fR}.$$ Also the
symmetric subgroup $N_{\bfb}:=Res_{k|\fQ}N_b'$ is defined over $\fQ$ and
satisfies $N_{\bfb}(\fR)=\prod_{\gs\in \gS_{\infty}}{^{\gs}(N_b')}_{\fR}$.
It follows that $(P_{\hbox{\scsi\bf b}},N_{\bfb})$ are incident: conditions
1) and 2) follow from the corresponding facts for $(P_b',N_b')$; we should
check 3). But since it is obvious that ${^{\gs}(}L_b')_{\fR}\inn
{^{\gs}(}P_b')_{\fR}$ is a hermitian Levi factor, the same holds for
$L_{\bfb}\inn P_{\bfb}$; 3) is satisfied.  This completes the proof of
\begin{proposition}\label{p16.2} Theorem \ref{t12.1} is true for the
  exceptional groups in the rank$\geq 2$, not split over $\fR$ case.
\end{proposition}

We are left with the classical cases. Here we may use the interpretation of
$G(\fR)$ as the unitary group of a $\pm$symmetric/hermitian form as in
(\ref{e7.1}), and $G$ is a $\fQ$-form of this. The precise realisation of
this is the interpretation in terms of central simple algebras with
involution; this is discussed in \cite{W}. More precisely, the algebraic
groups $G'$ which represent the indices of Proposition \ref{p16.1} are
(here we describe reductive groups; the corresponding derived groups are
the simple groups $G'$).
\begin{itemize}\item[(I)] \begin{itemize}\item[$D$:] degree $d$ central
    simple division algebra over $K$, $K|k$ an imaginary quadratic
    extension, $D$ has a $K|k$-involution (involution of the second kind).
  \item[$V$:] right $D$-vector space, of dimension $m$ over $D$, $dm=n+1$.
  \item[$h$:] hermitian form $h:V\times V\lra D$ of Witt index $s$, $2s\leq
    m$ ($2s<m$ if $d=1$), given by a matrix $H$.
  \item[$G'$:] unitary group $U(V,h)=\{g\in GL_D(V)|gHg^*=H\}$.
\end{itemize}
{\bf index:} ${^2A}^{(d)}_{n,s}$.
\item[(II)] \begin{itemize}\item[$D$:] totally definite quaternion division
    algebra, central simple over $k$, with canonical involution.
  \item[$V$:] right $D$-vector space of dimension $m$ over $D$.
  \item[$h$:] skew-hermitian form $h$ of Witt index $s< [{m\over 2}]$,
    given by a matrix $H$.
  \item[$G'$:] unitary group $U(V,h)=\{g\in GL_D(V)|gHg^*=H\}$.
\end{itemize}
{\bf index:} $D^{(2)}_{m,s}$ ($m$ even), ${^2D}^{(2)}_{m,s}$ ($m$ odd).
\item[(III)] \begin{itemize}\item[$D$:] totally indefinite quaternion
    division algebra, central simple over $k$, with the canonical
    involution.
  \item[$V$:] right $D$-vector space of dimension $m$.
  \item[$h$:] hermitian form $h:V\times V\lra D$ of Witt index $s$, $2s\leq
    m$, given by a matrix $H$.
  \item[$G'$:] unitary group $U(V,h)=\{g\in GL_D(V)|gHg^*=H\}$.
\end{itemize}
{\bf index:} $C^{(2)}_{m,s}$.
\end{itemize}

Finally, we must consider the following ``mixed cases'', which still can
give rise to groups of hermitian type:
\begin{itemize}\item[(II-IV):]\begin{itemize}\item[$D$:] a quaternion
    division algebra over $k$, with $D_{\nu}$ definite for
    $\nu_1,\ldots,\nu_a$, $D_{\nu}$ indefinite for
    $\nu_{a+1},\ldots,\nu_f$, where $f=[k:\fQ]$.
  \item[$V$:] same as for (II) above.
  \item[$h$:] same as for (II) above, $h$ of Witt index $s$.
  \item[$G'$:] same as for (II) above.
\end{itemize}\end{itemize}
$G(\fR)$ is then a product $(SU(n,\fH))^a\times (SO(2n-2,2))^{f-a}$, where
we have taken into account that we are assuming $G$ to be isotropic and of
hermitian type.  Note however, that since the factors $SO(2n-2,2)$
corresponding to the primes $\nu_{a+1},\ldots,\nu_f$ have $\fR$-split torus
of dimension two, the $k$-rank of $G'$ must be $\leq2$. Hence the only
indices where this can occur are: ${^iD}^{(2)}_{n,1}$ and
${^i}D^{(2)}_{n,2}$, $i=1,2$.

In terms of the spaces $(V,h)$, the standard parabolics are stabilizers of
totally isotropic subspaces $H_b\inn V$, where $H_1$ is one-dimensional
(over $D$), while $H_s$ is a maximal totally isotropic subspace. The latter
case corresponds to zero-dimensional boundary components. We consider first
the case $H_b,\ b<s$, of which at least $H_1$ exists, because of the
assumption rank $\geq 2$. Fix a basis $h_1,\ldots, h_b$ of $H_b$ of
isotropic vectors $h(h_i,h_i)=0$ for all $i=1,\ldots,b$. Then there exist,
in $V$, elements $h_i',\ i=1,\ldots,b$ with $h(h_i,h_j')=\gd_{ij}$, and
$h_1',\ldots,h_b'$ span a complementary totally isotropic subspace; denote
it by $H_b'$. Then $H:=H_b\oplus H_b'$ is a {\it non-degenerate} space for
$h$, $h_{|H}$ is a non-degenerate form. It follows that $\{g\in GL(V) |
g(H)\inn H\} = \{g\in GL(V) | g(H^{\perp})\inn H^{\perp}\}$. In the
following we will work in the (reductive) unitary group $G'=U(V,h)$; for
any subgroup $H\inn G'$ we can take the intersection $SL(V)\cap H\inn
SL(V)\cap G'$ to give subgroups of the simple group. Furthermore, up to
Corollary \ref{c19.1} below, we omit the primes in the notations for the
subgroups of $G'$. Set
\begin{equation}\label{e18.1} N=U(H,H^{\perp};h)=\{g\in GL(V)|g(H)\inn H,\
  g(H^{\perp})\inn H^{\perp}\};
\end{equation}
then $N=U(H,h_{|H})\times U(H^{\perp},h_{|H^{\perp}}),$ and
$U(H,h_{|H})\cong \cZ_G(U(H^{\perp},h_{|H^{\perp}}))$. So setting
$L=U(H^{\perp}, h_{|H^{\perp}})$, we have
\begin{equation}\label{e18.2} N\cong L \times \cZ_G(L).
\end{equation}

Next we note that the basis $h_1,\ldots, h_b$ of $H_b$ determines a unique
$\fR$-split torus $A_b\inn A$, where $A$ is the maximal $\fR$-split torus
defined by a basis $h_1,\ldots,h_s$ of a maximal totally isotropic subspace
$H_s\nni H_b$, namely the scalars $\ga=\ga\cdot{\bf 1}\in GL(H_b)$,
extended to $GL(V)$ by unity. Taking the centralizer of the torus $A_b$
gives a Levi factor of the parabolic $P_b=\cN_G(H_b),\ b<s$ (the normalizer
in $G$ of $H_b$).
\begin{lemma}\label{l18.1} $L=U(H^{\perp},h_{|H^{\perp}})$ is the
  $k$-hermitian factor $G_b^{(1)}=M_b\cdot L_b$ of $P_b$ in the
  decomposition of $P_b$ as in Theorem \ref{t4.1}.
\end{lemma}
{\bf Proof:} First observe that $L\inn P_b$, as $H^{\perp}$ is orthogonal
to the totally isotropic subspace, hence $L$ normalizes $H_b$. Since $L$ is
reductive, there is a Levi decomposition of $P_b$ for which $L$ is
contained in the Levi factor. It is clearly of hermitian type, and maximal
with this property. We must explain why the Levi factor is the standard one
$\cZ(A_b)$. But this follows from the fact that $H_b$ is constructed by
means of a basis, which in turn was determined by the choice of $\fR$-split
torus $A_b$. It therefore suffices to explain the ``compact'' factor $M_b$.
This factor occurs only in the cases $\bf I_{\bf\scst p,q}$ and $\bf
IV_{\scst\bf n}$. We don't have to consider the latter case, as this is
split over $\fR$ if rank $\geq 2$.  So suppose $G\cong U(V,h)$, where
$(V,h)$ is as in (I) above. We first determine the anisotropic kernel. Let
$H_s$ be a maximal totally isotropic subspace, $S:=H_s\oplus H_s'$ as
above. Then $U(S^{\perp}, h_{|S^{\perp}})$ is the anisotropic kernel,
$U(S^{\perp},h_{|S^{\perp}})(\fR)\cong U(md-2sd)$.  In particular, for
$m=2s$, there is no anisotropic kernel. Now consider the group
$L=U(H^{\perp},h_{|H^{\perp}})$. Clearly, for $b<s$, we have
\[U(S^{\perp},h_{|S^{\perp}})\inn U(H^{\perp},h_{|H^{\perp}})=L,\]
so that $L$ contains the anisotropic kernel. Note that
$SU(H^{\perp},h_{|H^{\perp}})(\fR)\cong L_b(\fR)$, while (if $H^{\perp}\neq
\{0\}$)
\[U(H^{\perp},h_{|H^{\perp}})(\fR)/SU(H^{\perp},h_{|H^{\perp}})(\fR)\cong
M_b(\fR)\cong U(1).\] Here we have used that
$U(H^{\perp},h_{|H^{\perp}})\inn SU(V,h)$, as it is for the group $SU(V,h)$
(and not for $U(V,h)$) that $M_b(\fR)\cong U(1)$ (see \cite{S}, p.~115).
This verifies the Lemma for the groups of type $\bf I$. \ende

Now note that $L(\fR)\cong(M_b\cdot L_b)(\fR)=M_b(\fR)L_b(\fR)$, so for the
domain defined by $L$ we have
$\cD_L=M_b(\fR)L_b(\fR)/M_b(\fR)K_b=L_b(\fR)/K_b\cong F_b$, hence
$\cD_N\cong \cD_{N_b}$ as in Table \ref{T1}. Consider also $\cZ_G(L_b)$ and
$\cZ_G(M_bL_b)$; both are defined over $\fR$, and clearly
$\cZ_{G(\fR)}(L_b(\fR))/M_b(\fR)\cong \cZ_{G(\fR)}(M_b(\fR)L_b(\fR))$, so
the group $L\times \cZ_G(L)$ (both these factors being defined over $k$)
is, over $\fR$,
\begin{eqnarray}\label{e19.1}
  L(\fR)\times \cZ_G(L)(\fR) & \cong & M_b(\fR)\cdot L_b(\fR)\times
  \cZ_{G(\fR)}(M_b(\fR)L_b(\fR)) \\ & \cong & M_b(\fR)\cdot L_b(\fR)\times
  \cZ_{G(\fR)}(L_b(\fR))/M_b(\fR) \nonumber \\ & \cong & L_b(\fR)\times
  \cZ_{G(\fR)}(L_b(\fR)).\nonumber
\end{eqnarray}
This completes the proof of
\begin{proposition}\label{p19.1} The subgroup $N$ of (\ref{e18.1})
  satisfies $N(\fR)\cong N_b(\fR)$, the latter group being the standard
  symmetric subgroup (\ref{e10.1}) standard incident to $P_b(\fR)$.
  \end{proposition}
\begin{corollary}\label{c19.1} The parabolic $P_b$ and the symmetric
  subgroup $N$ of (\ref{e18.1}) are incident over $k$, i.e.,
  $(P_b(\fR),N(\fR))$ are incident in the sense of Definition \ref{d9.1}.
\end{corollary}
Up to this point we have been working with the absolutely simple $k$-group;
we now denote this situation by $G'$ as in section 3.1, and consider
$G=Res_{k|\fQ}G'$. Let again primes in the notations denote subgroups of
$G'$, the unprimed notations for subgroups of $G$. As above we set
$P_{\bfb}:= Res_{k|\fQ}(P_b')$, and we denote the subgroup $N'$ of
(\ref{e18.1}) henceforth by $N_b'$ and set: $N_{\bfb}:=Res_{k|\fQ}(N_b')$.
Then Corollary \ref{c19.1} tells us that $(P_b'(\fR),N'(\fR))$ are
incident. We now claim
\begin{lemma}\label{l19.1} The $\fQ$-groups $(P_{\bfb},N_{\bfb})$ are incident.
\end{lemma}
{\bf Proof:} $P_{\bfb}(\fR)$ is a product $P_{\bfb}(\fR)\cong
P_{b,1}(\fR)\times \cdots \times P_{b,d}(\fR)$ corresponding to
(\ref{e9.1}); by assumption $F_{\bfb}$ is not zero-dimensional. Hence for
at least one factor $P_{b,\gs}(\fR)$ the incident group
${^{\gs}N}'_b(\fR)\cong{^{\gs}L}'_b(\fR)\times
{\cZ}_{{^{\gs}G}'(\fR)}({^{\gs}L}'_b(\fR))$ is defined. Consequently
$N_{\bfb}$ is not trivial, and it is clearly a $\fQ$-form of
$N_{\bfb}(\fR)=\prod_{\gs}{^{\gs}N}'_b(\fR)$. \ende With Corollary
\ref{c19.1} and Lemma \ref{l19.1}, we have just completed the proof of the
following.
\begin{proposition}\label{p19.2} To each standard maximal $\fQ$-parabolic
  $P_{\bfb}$ of $G$ with $b<s$, there is a symmetric $\fQ$-subgroup
  $N_{\bfb}\inn G$ such that $(P_{\bfb},N_{\bfb})$ are incident.
\end{proposition}

Finally, we remark on what happens for the parabolic corresponding to the
zero-dimensional boundary components. We have, in the notations above,
$H=H_s\oplus H_s'$, and $H^{\perp}$ is anisotropic for $h$.  It follows
that the group $L$ of Lemma \ref{l18.1} is anisotropic; its semisimple part
is the semisimple anisotropic kernel of $G'$. If $s={1\over 2}\dim V$, then
$H=V$ already, $H^{\perp}=\{0\}$. Hence the group $N$ of (\ref{e18.2}) is
the whole group ($L=1 \Ra \cZ_G(L)=G$). Otherwise it is of the form
$\{\hbox{anisotropic}\}\times \{\hbox{$k$-split}\}$. We list these in Table
\ref{T3}. Note that the domains occuring have $\fR$-rank equal to the
$\fQ$-rank of $G$, suggesting this as a possible modification of the
definition of incident:
\begin{itemize}\item[1')] $N$ has $\fR$-rank equal to the $\fQ$-rank of
  $G$.
\end{itemize}
Viewing things this way, we see that again indices $C^{(1)}$ represent an
exception; for these 1) and 1') are equivalent.

\begin{table}\caption{\label{T3} $k$-subgroups incident with
    zero-dimensional boundary components}
  $$\begin{array}{|c|c|c|c|c|} \hline \hbox{Index} & L & \cZ_G(L) &
    \hbox{subdomains} & \cZ_G(L)(\fR) \\ \hline \hline {^2A}^{(d)}_{n,s} &
    {^2A}^{(d)}_{n-2ds,0} & {^2A}^{(d)}_{2ds-1,s} & \bf
    I_{\hbf{p-ds,q-ds}}\times I_{\hbf{ds,ds}} & SU(ds,ds) \\ \hline
    {^1D}^{(2)}_{n,s} & {^1D}^{(2)}_{n-2s,0} & {^1D}^{(2)}_{2s,s} & \bf
    II_{\hbf{n-s}}\times II_{\hbf{s}} & SU(2s,\fH)\ (\hbox{$n$ even}) \\
    \hline {^2D}^{(2)}_{n,s} & {^2D}^{(2)}_{n-2s,0} & {^2D}^{(2)}_{2s,s} &
    \bf II_{\hbf{n-s}}\times II_{\hbf{ s}} & SU(2s,\fH)\ (\hbox{$n$ odd})
    \\ \hline C^{(1)} & - & G & - & - \\ \hline C^{(2)}_{n,s} &
    C^{(2)}_{n-s,0} & C^{(2)}_{s,s} & \bf III_{\hbf{n-s}}\times
    III_{\hbox{\scsi\bf s}} & Sp(2s,\fR) \\ \hline
\end{array}$$
\end{table}

Let us now see which of the subgroups listed in Table \ref{T2} are defined
over $k$. We use the notations $D,\ V, h$ and $G$ as described above in the
cases (I)-(III).
\begin{itemize}\item[(I)] Again $d$ denotes the degree of $D$. In $U(V,h)$
  we have the subgroup $U(V',h_{|V'})$ for any codimension one subspace
  $V'\inn V$. Let $W=(V')^{\perp}$ be the one-dimensional (over $D$)
  subspace orthogonal to $V'$. Then $U(W,h_{|W})$ is again a unitary group
  whose set of $\fR$-points is isomorphic to $U(p_W,q_W)$ for some
  $p_W,q_W$. Actually each $h_{\nu}$ for each infinite prime $\nu$ gives an
  $\fR$-group $U(p_{\scstst W,\nu},q_{\scstst W,\nu})$.  Let
  $(p_{\nu},q_{\nu})$ be the signature of $h_{\nu}$ on $V_{\nu}$. Then
  $U(V'_{\nu},{h_{\nu}}_{|V'_{\nu}})\cong U(p_{\nu}-p_{\scstst
    W,\nu},q_{\nu}-q_{\scstst W,\nu})$.  This gives rise to a product
  $N=\prod U(p_{\scstst W,\nu},q_{\scstst W,\nu})\times
  U(p_{\nu}-p_{\scstst W,\nu},q_{\nu}-q_{\scstst W,\nu})$, and the factors
  of the domain $\cD_N$ are of type $\bf I_{\scst\bf p_{\hbox{$\scstst
        W$},\nu},q_ {\hbox{$\scstst W$},\nu}}\times I_{\scst\bf
    p_{\nu}-p_{\hbox{$\scstst W$},\nu},q_{\nu}-q_{\hbox{$\scstst
        W$},\nu}}$.  In particular, for $p_{\scstst W,\nu}=0$, this is an
  irreducible group of type $\bf I_{\scst\bf
    p_{\nu},q_{\nu}-q_{\hbox{$\scstst W$},\nu}}$ and for $q_{\scstst
    W,\nu}=0$, of type ${\bf I}_{\scst\bf p_{\nu}-p_{\hbox{$\scstst
        W$},\nu},q_{\nu}}$.  Now since $k$ is the degree of $D$, all of
  $p_{\nu}, q_{\nu}, p_{\hbox{$\scstst W$},\nu}, q_{\hbox{$\scstst
      W$},\nu}$ are divisible by $d$ and the net subdomains these subgroups
  (possibly) define are
\begin{equation}\label{eZZ}
  \bf I_{\hbf{p-jd,q}},\quad I_{\hbf{p,q-jd}},\quad I_{\hbf{p-id,q-jd}}
  \inn I_{\hbf{p,q}}, \ \ i,j=1,\ldots, s.
\end{equation}

\item[(II)] In $U(V,h)$ we have as above $U(V',h_{|V'})$; now if $h$ is
  non-degenerate on $V'$, then $U(V',h_{|V'})\cong U(n-1,D)$, giving
  subgroups of the real groups, defined over $k$, of type $U(n-1,\fH)\inn
  U(n,\fH)$, with a corresponding subdomain of type $\bf II_{\scst\bf
    n-1}\inn II_{\scst\bf n}$. This occurs at the primes for which $D$ is
  definite; at the others $SU(V',h_{|V'})\inn SU(V,h)$ is of the type
  $SO(2n-4,2)\inn SO(2n-2,2)$ (for $n$=dimension of $V$ over $D$). So we
  have maximal $k$-domains
  $$\bf II_{\hbf{n-1}}\inn II_{\hbf{n}},\ (\nu \hbox{ definite}),\quad
  \quad IV_{\hbf{2n-4}}\inn IV_{\hbf{2n-2}},\ (\nu \hbox{ indefinite}).$$
\item[(III)] The index is $C^{(2)}_{n,s}$; this case in considered in more
  detail below.
\end{itemize}
{}From this, we deduce
\begin{proposition}\label{p23} Let $G'$ have $\rank_kG'=s\geq2$, not split
  over $\fR$, and let $P_s'$ be a standard $k$-parabolic defining a
  zero-dimensional boundary component, $P_s'(\fR)=N(F)$, and $\dim(F)=0$.
  Then there is a $k$-subgroup $N'$ incident with $P_s'$, with the
  following exception: Index $C^{(2)}_{2s,s}$.
\end{proposition}
{\bf Proof:} We first deduce for which of the indices listed in Proposition
\ref{p16.1} zero-dimensional boundary components of $\cD'$ are rational
(this is necessary for the zero-dimensional boundary components of $\cD$ to
be rational). We need not consider exceptional cases or type $\bf
IV_{\bfn}$.  We first consider the groups of type ${^2A}$.
\begin{Lemma}\label{L19A} For $G'$ with the index ${^2A}^{(d)}_{n,s}$, let
  $G'(\fR)\cong SU(p,q)$. Then the zero-dimensional boundary components are
  rational $\iff$ $sd=q$.
\end{Lemma}
{\bf Proof:} Let $H_s$ be an $s$-dimensional (over $D$) totally isotropic
subspace, with basis $h_1,\ldots,h_s$. Let $h_i'\in V$ be vectors such that
$h(h_i,h_j')=\gd_{ij}$, and set $H_s'=<h_1',\ldots,h_s'>$. Then $h$,
restricted to $H:=H_s\oplus H_s'$, is non-degenerate, and
$SU(H^{\perp},h_{|H^{\perp}})$ is the anisotropic kernel. The group
$SU(H,h_{|H})(\fR) \cong SU(sd,sd)$, while
$SU(H^{\perp},h_{|H^{\perp}})(\fR)\cong SU(p-ds,q-ds)$. This defines the
subdomain of type $\bf I_{\hbf{ds,ds}}\times I_{\hbf{p-ds,q-ds}}$ of Table
\ref{T1}, hence the boundary component, which is the second factor, is
zero-dimensional $\iff$ $q=ds$. \ende As to indices of type $D$ we observe
the following.
\begin{Lemma}\label{l5.8.1} $\dim(F)=0$ does not occur for the indices
  of type (II) in Proposition \ref{p16.1}.
\end{Lemma}
{\bf Proof:} Recall that $D$ is a quaternion division algebra, central
simple over $k$, with the canonical involution, $V$ is an $n$-dimensional
right $D$-vector space, and $h:V\times V\lra D$ is a skew-hermitian form.
Let $\nu_1,\ldots, \nu_a$ denote the infinite primes for which $D_{\nu}$ is
definite, $\nu_{a+1},\ldots, \nu_d$ the primes at which $D_{\nu}$ is split.
Then $G(\fR)$ is a product
$$(SU(n,\fH))^a\times (SO(2n-2,2))^{d-a},$$ where we have taken into
account that $G$ is assumed to be of hermitian type. At each of the first
factors we have the Satake diagram
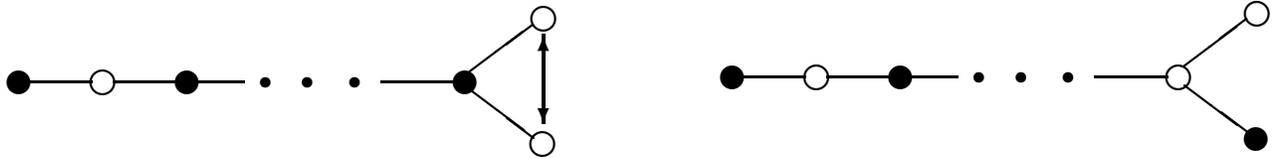
\setlength{\unitlength}{0.005500in}%
$$\begin{picture}(1202,145)(69,691) \thicklines \put(240,760){\circle*{22}}
  \put(160,760){\circle{22}} \put(315,760){\circle*{10}}
  \put(355,760){\circle*{10}} \put(400,760){\circle*{10}}
  \put(505,760){\circle*{22}} \put(579,701){\circle{22}} \put(
  90,760){\line( 1, 0){ 60}} \put(170,760){\line( 1, 0){ 60}}
  \put(250,760){\line( 1, 0){ 45}} \put(425,760){\line( 1, 0){ 70}}
  \put(510,770){\line( 4, 3){ 60}} \put(511,752){\line( 4,-3){ 60}}
  \put(580,805){\vector( 0,-1){ 85}} \put(580,805){\vector( 0,1){ 0}} \put(
  80,760){\circle*{22}} \put(1191,757){\line( 4,-3){ 60}}
  \put(580,820){\circle{22}} \put(1260,825){\circle{22}}
  \put(760,765){\circle*{22}} \put(920,765){\circle*{22}}
  \put(840,765){\circle{22}} \put(995,765){\circle*{10}}
  \put(1035,765){\circle*{10}} \put(1080,765){\circle*{10}}
  \put(1185,765){\circle{22}} \put(1259,706){\circle*{22}}
  \put(770,765){\line( 1, 0){ 60}} \put(850,765){\line( 1, 0){ 60}}
  \put(930,765){\line( 1, 0){ 45}} \put(1105,765){\line( 1, 0){ 70}}
  \put(1190,775){\line( 4, 3){ 60}}
\end{picture}$$

\hspace*{3cm} for $n$ odd, \hspace*{6.5cm} for $n$ even.

The corresponding $\fR$-root systems are then:

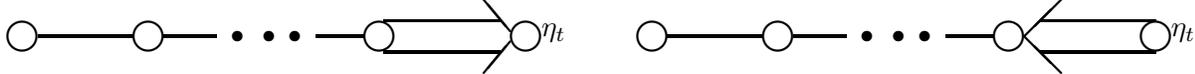
\setlength{\unitlength}{0.005500in}%
\begin{picture}(1108,70)(16,700)
  \thicklines \put(150,735){\circle{28}} \put(370,735){\circle{28}}
  \put(510,735){\circle{28}} \put(525,735){$\eta_t$}
  \put(235,735){\circle*{10}} \put(265,735){\circle*{10}}
  \put(290,735){\circle*{10}} \put( 45,735){\line( 1, 0){ 90}}
  \put(165,735){\line( 1, 0){ 50}} \put(310,735){\line( 1, 0){ 45}}
  \put(470,770){\line( 5,-6){ 25}} \put(495,730){\line(-5,-6){ 25}}
  \put(375,720){\line( 1, 0){110}} \put(375,750){\line( 1, 0){110}} \put(
  30,735){\circle{28}} \put(630,735){\circle{28}} \put(1000,750){\line( 1,
    0){110}} \put(750,735){\circle{28}} \put(970,735){\circle{28}}
  \put(1110,735){\circle{28}} \put(1125,735){$\eta_t$}
  \put(835,735){\circle*{10}} \put(865,735){\circle*{10}}
  \put(890,735){\circle*{10}} \put(645,735){\line( 1, 0){ 90}}
  \put(765,735){\line( 1, 0){ 50}} \put(910,735){\line( 1, 0){ 45}}
  \put(1000,720){\line( 1, 0){110}} \put(985,735){\line( 1,-1){ 35}}
  \put(985,735){\line( 1, 1){ 35}}
\end{picture}

In particular, the $\fR$-root corresponding to the parabolic $P_{t}$ with
$\dim(F_{t})=0$ is the right-most one. On the other hand, the $k$-index is
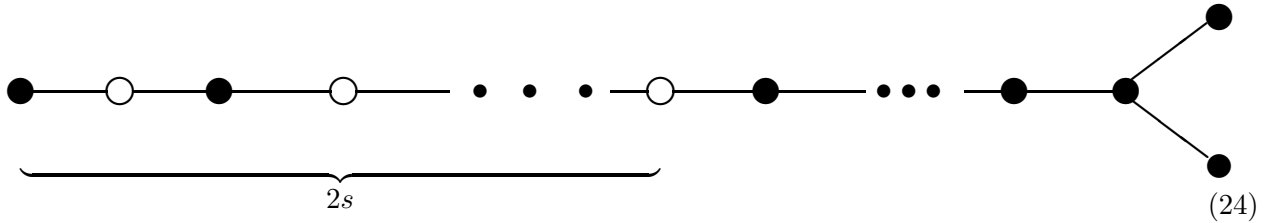
\begin{equation}\label{eZ1.1}
\setlength{\unitlength}{0.006500in}%
\begin{picture}(987,141)(14,670)
  \thicklines \put(25,680){$\underbrace{\hspace*{8.5cm}}_{\hbox{$2s$}}$}
  \put(435,740){\circle*{10}} \put(480,740){\circle*{10}} \put(
  25,740){\circle*{22}} \put(185,740){\circle*{22}}
  \put(105,740){\circle{22}} \put(285,740){\circle{22}}
  \put(540,740){\circle{22}} \put(625,740){\circle*{22}}
  \put(740,740){\circle*{10}} \put(760,740){\circle*{10}}
  \put(720,740){\circle*{10}} \put(825,740){\circle*{22}}
  \put(395,740){\circle*{10}} \put(915,740){\circle*{22}}
  \put(921,732){\line( 4,-3){ 60}} \put(990,680){\circle*{20}}
  \put(990,800){\circle*{22}} \put( 35,740){\line( 1, 0){ 60}}
  \put(115,740){\line( 1, 0){ 60}} \put(195,740){\line( 1, 0){ 80}}
  \put(295,740){\line( 1, 0){ 75}} \put(500,740){\line( 1, 0){ 30}}
  \put(550,740){\line( 1, 0){ 65}} \put(635,740){\line( 1, 0){ 70}}
  \put(705,740){\line(-1, 0){ 5}} \put(785,740){\line( 1, 0){ 40}}
  \put(835,740){\line( 1, 0){ 70}} \put(920,750){\line( 4, 3){ 60}}
\end{picture}
\end{equation}
with the $k$-root system

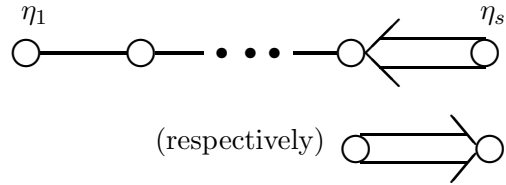
\setlength{\unitlength}{0.00500in}%
$$\begin{picture}(513,170)(51,610) \thicklines \put(185,745){\circle{28}}
  \put(405,745){\circle{28}} \put(545,745){\circle{28}}
  \put(540,780){$\eta_s$} \put(270,745){\circle*{10}}
  \put(300,745){\circle*{10}} \put(325,745){\circle*{10}} \put(
  80,745){\line( 1, 0){ 90}} \put(200,745){\line( 1, 0){ 50}}
  \put(345,745){\line( 1, 0){ 45}} \put( 65,745){\circle{28}} \put(
  60,780){$\eta_1$} \put(435,730){\line( 1, 0){110}} \put(415,660){\line(
    1, 0){110}} \put(420,745){\line( 1,-1){ 35}} \put(420,745){\line( 1,
    1){ 35}} \put(435,760){\line( 1, 0){110}} \put(200,645){(respectively)}
  \put(410,645){\circle{28}} \put(550,645){\circle{28}}
  \put(510,680){\line( 5,-6){ 25}} \put(535,640){\line(-5,-6){ 25}}
  \put(415,630){\line( 1, 0){110}}
\end{picture}$$
from which it is evident that $P_{t}$ is defined over $k$ $\iff$ $s=t$
($=[{n\over 2}]$).  But this is the split over $\fR$ case. Consequently,
$a=0$ and $D$ is totally indefinite.

So we consider a prime $\nu$ where $D_{\nu}$ is split; the $\fR$-index is
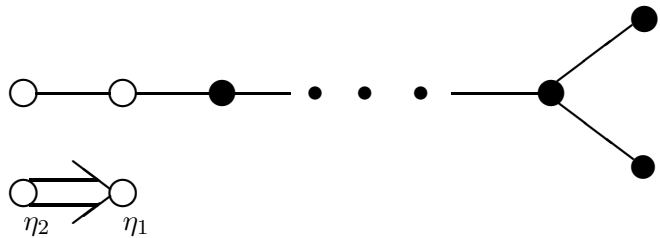
\setlength{\unitlength}{0.006500in}%
$$\begin{picture}(522,171)(69,660) \thicklines \put(315,760){\circle*{10}}
  \put(355,760){\circle*{10}} \put(400,760){\circle*{10}}
  \put(505,760){\circle*{22}} \put(580,820){\circle*{22}}
  \put(579,701){\circle*{20}} \put( 80,760){\circle{22}}
  \put(160,760){\circle{22}} \put( 80,680){\circle{22}}
  \put(240,760){\circle*{22}} \put(160,680){\circle{22}}
  \put(120,655){\line( 4, 3){ 30}} \put(80,650){$\eta_2$} \put(
  90,760){\line( 1, 0){ 60}} \put(170,760){\line( 1, 0){ 60}}
  \put(250,760){\line( 1, 0){ 45}} \put(425,760){\line( 1, 0){ 70}}
  \put(510,770){\line( 4, 3){ 60}} \put(511,752){\line( 4,-3){ 60}} \put(
  85,690){\line( 1, 0){ 55}} \put( 85,670){\line( 1, 0){ 55}}
  \put(120,705){\line( 4,-3){ 30}} \put(160,650){$\eta_1$}
\end{picture}
$$ the $\fR$-root $\eta_2$ corresponding to the two-dimensional totally
isotropic subspace and zero-dimensional boundary component. The $k$-index
is as in (\ref{eZ1.1}), so $\eta_2$ is always anisotropic; the boundary
components are actually one-dimensional.  This verifies the statements of
the lemma. \ende Note that this proves Proposition \ref{p23} for the
indices of type (II).

Now consider index $C^{(2)}_{n,s}$. The $k$-index is
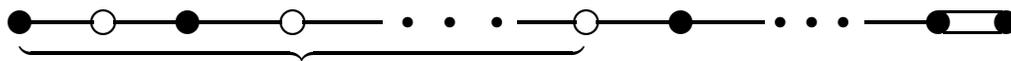
\setlength{\unitlength}{0.005500in}%
$$\begin{picture}(962,22)(14,729) \thicklines \put(435,740){\circle*{10}}
  \put(480,740){\circle*{10}} \put( 25,740){\circle*{22}} \put(
  25,720){$\underbrace{\hspace*{7.5cm}}_{\hbox{$2s$}}$}
  \put(185,740){\circle*{22}} \put(105,740){\circle{22}}
  \put(285,740){\circle{22}} \put(565,740){\circle{22}}
  \put(655,740){\circle*{22}} \put(750,740){\circle*{10}}
  \put(780,740){\circle*{10}} \put(810,740){\circle*{10}}
  \put(395,740){\circle*{10}} \put(900,740){\circle*{22}}
  \put(901,750){\line( 1, 0){ 60}} \put(965,740){\circle*{22}} \put(
  35,740){\line( 1, 0){ 60}} \put(115,740){\line( 1, 0){ 60}}
  \put(195,740){\line( 1, 0){ 80}} \put(295,740){\line( 1, 0){ 75}}
  \put(500,740){\line( 1, 0){ 55}} \put(575,740){\line( 1, 0){ 70}}
  \put(665,740){\line( 1, 0){ 70}} \put(735,740){\line(-1, 0){ 5}}
  \put(830,740){\line( 1, 0){ 60}} \put(905,730){\line( 1, 0){ 55}}
\end{picture}$$
and the $k$-root system is
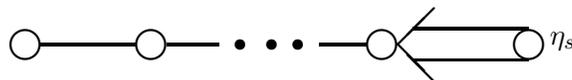
\setlength{\unitlength}{0.005500in}%
$$\begin{picture}(508,70)(51,710) \thicklines \put(185,745){\circle{28}}
  \put(405,745){\circle{28}} \put(545,745){\circle{28}}
  \put(565,745){$\eta_s$} \put(270,745){\circle*{10}}
  \put(300,745){\circle*{10}} \put(325,745){\circle*{10}} \put(
  65,745){\circle{28}} \put( 80,745){\line( 1, 0){ 90}}
  \put(435,760){\line( 1, 0){110}} \put(200,745){\line( 1, 0){ 50}}
  \put(345,745){\line( 1, 0){ 45}} \put(435,730){\line( 1, 0){110}}
  \put(420,745){\line( 1,-1){ 35}} \put(420,745){\line( 1, 1){ 35}}
\end{picture}$$
The same reasoning as above shows that $F_{t}$ is rational $\iff$ $2s=t$,
but that is only possible if the index is $C^{(2)}_{2n,n}$. Hence:
\begin{Lemma}\label{l5.8.2}
  The only indices of Proposition \ref{p16.1}, case (III), for which
  zero-dimensional boundary components occur are $C^{(2)}_{2n,n}$.
\end{Lemma}
This index is that of the unitary group $U(V,h)$, where $V$ is a
$2n$-dimensional vector space over $D$, and $h$ has Witt index $n$. We can
find $n$ hyperbolic planes $V_i$ such that
\[V=V_1\oplus\cdots \oplus V_n.\]
This decomposition is defined over $k$, hence the subgroup
\[N=U(V_1,h_{|V_1})\times \cdots \times U(V_n,h_{|V_n}),\]
which is a product of groups with index $C^{(2)}_{2,1}$, is also defined
over $k$. We have
\begin{equation}\label{E20}N(\fR)\cong \underbrace{Sp(4,\fR)\times \cdots
    \times Sp(4,\fR)}_{n\ \hbox{\scsi times}}
\end{equation}
and the domain $\cD_N$ is of type $\bf (III_{\hbox{\scsi\bf 2}})^n$. This
is the exception in the statement of the main theorem.

\vspace*{.3cm}
\noindent{\bf Proof of Proposition \ref{p23}:} We have already completed
the proof for (II) and (III), and as we mentioned above, the exceptional
cases and (IV) need not be considered. It remains to show the existence of
groups of the stated types for indices ${^2A}$. We explained above how one
can find $k$-subgroups $N$ such that $\cD_N$ has irreducible components of
types $\bf I_{\hbox{\scsi\bf p-jd,q}}$ (see (\ref{eZZ})). Here we take a
maximal totally isotropic subspace $H_s$, and $H:=H_s\oplus H_s'$ as
described there. Let $H^{\perp}$ denote the orthogonal complement, so that
$SU(H^{\perp},h_{|H^{\perp}})$ is the anisotropic kernel. Then, if
$G'(\fR)=SU(p,q)$, we have
\[SU(H,h_{|H})(\fR)\cong SU(sd,sd),\quad SU(H^{\perp},h_{|H^{\perp}})\cong
SU(p-sd,q-sd).\] Therefore we get a subdomain of type \[\bf
I_{\hbf{sd,sd}}\times I_{\hbf{p-sd,q-sd}},\] which is irreducible $\iff$
$sd=q$; Then $N=\{g\in G\Big| g(H)\subseteq H\}$ is a $k$-subgroup with
$N(\fR)\sim SU(q,q)\times\{\hbox{compact}\}$, and $N$ then fulfills 1), 2')
and 3'). By Lemma \ref{L19A}, this holds precisely when the boundary
component $F_s$ is a point. This completes the proof if $p>q$. It remains
to consider the case where $\cD'$ is of type $\bf I_{\hbf{q,q}}$. In this
case, $q=d\cdot j$ for some $j$, and the hermitian form $h:V\times V\lra D$
has Witt index $j$. The vector space $V$ is then $2j$-dimensional, and it
is the orthogonal direct sum of hyperbolic planes, $V=V_1\oplus \cdots
\oplus V_j$, $\dim_DV_i=2$. Consider the $k$-subgroup
\[ N=\{g\in GL_D(V) \Big| g(V_i)\inn V_i, i=1,\ldots,j\}.\]
Clearly $N\cong N_1\times \cdots \times N_j$, and each $N_i$ is a subgroup
of rank one with index ${^2A}^{(d)}_{2d-1,1}$. As was shown in \cite{hyp},
in each $N_i$ we have a $k$-subgroup $N_i'\inn N_i$, with $\cD_{N_i'}$ of
type $({\bf I_{\hbf{1,1}}})^d$. Then
\[ N':=N_1'\times \cdots \times N_j'\]
is a $k$-subgroup with $\cD_{N'}$ of type $(({\bf I_{\hbf{1,1}}})^d)^j=
({\bf I_{\hbf{1,1}}})^{d\cdot j}= ({\bf I_{\hbf{1,1}}})^q$, which is a
maximal polydisc, i.e., satisfies 1), 2'') and 3''). This completes the
proof of Proposition \ref{p23} in this case also. \ende

\section{Rank one}
We now come to the most interesting and challenging case. In this last
paragraph $G'$ will denote an absolutely simple $k$-group, $G$ the
corresponding $\fQ$-simple group, both assumed to have rank one. There is
only one standard maximal parabolic $P_1'\inn G'$ in this case, so we may
delete the subscript $_1$ in the notations. Let $P\inn G$ be the
corresponding $\fQ$-parabolic, so $P(\fR)=P_1(\fR)\times \cdots \times
P_d(\fR)$, where $P_{\nu}(\fR)\inn {^{\gs_{\nu}}G}'(\fR)$ is a standard
maximal parabolic, say $P_{\nu}(\fR)=N(F_{b_{\nu}}),\ F_{b_{\nu}}\inn
\overline{\cD}_{\gs_{\nu}}$.  As we observed above, the $F_{b_{\nu}}$ are
all hermitian spaces whose automorphism group is an $\fR$-form of some
fixed algebraic group. As we are now assuming the rank to be one, it
follows from Lemma \ref{L12a} that $L$ (=$L_{\bfb}$ in the notations above)
is anisotropic. One way that this may occur was explained there, namely
that if one of the factors $F_{b_{\nu}}$ is a {\it point}, in which case
the symmetric space of $L(\fR)$ has a compact factor. Another possibility
is that all $F_{b_{\nu}}$ are positive-dimensional, in which case $L$ is a
``genuine'' anisotropic group. The type of $F_{b_{\nu}}$ can be determined
from the $k$-index of $G'$ and the $\fR$-index of ${^{\gs}_{\nu}G}'$.  For
example, for $G'$ of type ${^2A}$, these indices are:

\vspace*{.5cm}
\setlength{\unitlength}{0.004500in}%
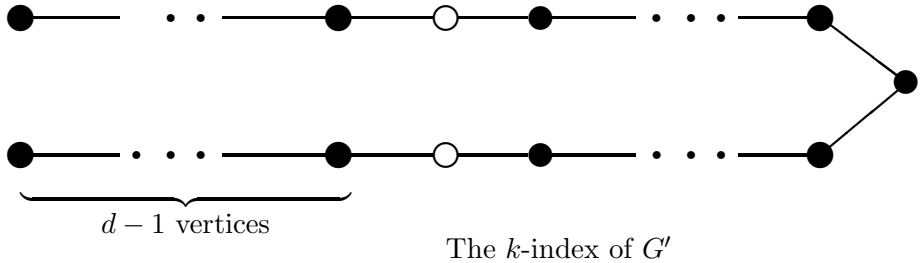
\begin{picture}(600,560)(90,235)
  \thicklines \put(280,780){\circle*{10}} \put(315,780){\circle*{10}}
  \put(105,780){\circle*{30}} \put(475,780){\circle*{30}}
  \put(120,780){\line( 1, 0){100}} \put(340,780){\line( 1, 0){120}}
  \put(845,780){\circle*{10}} \put(885,780){\circle*{10}}
  \put(920,780){\circle*{10}} \put(600,780){\circle{28}}
  \put(710,780){\circle*{28}} \put(1035,780){\circle*{30}}
  \put(490,780){\line( 1, 0){ 95}} \put(615,780){\line( 1, 0){ 80}}
  \put(720,780){\line( 1, 0){100}} \put(940,780){\line( 1, 0){100}}
  \put(240,620){\circle*{10}} \put(280,620){\circle*{10}}
  \put(315,620){\circle*{10}} \put(105,620){\circle*{30}}
  \put(475,620){\circle*{30}} \put(120,620){\line( 1, 0){100}}
  \put(340,620){\line( 1, 0){120}} \put(845,620){\circle*{10}}
  \put(885,620){\circle*{10}} \put(920,620){\circle*{10}}
  \put(600,620){\circle{28}} \put(710,620){\circle*{28}}
  \put(1035,620){\circle*{30}} \put(490,620){\line( 1, 0){ 95}}
  \put(615,620){\line( 1, 0){ 80}} \put(720,620){\line( 1, 0){100}}
  \put(940,620){\line( 1, 0){100}} \put(1135,705){\circle*{28}}
  \put(1035,780){\line( 4,-3){100}} \put(1035,620){\line( 6, 5){ 90}}

  \put(105,580){$\underbrace{\hspace*{4.4cm}}_{ \hbox{$d-1$ vertices}}$}
  \put(600,500){The $k$-index of $G'$}
\end{picture}

\setlength{\unitlength}{0.004500in}%
\begin{picture}(1400,60)(90,235)
  \thicklines

  \put(240,250){\circle*{10}} \put(280,250){\circle*{10}}
  \put(315,250){\circle*{10}} \put(840,410){\circle*{10}}
  \put(880,410){\circle*{10}} \put(915,410){\circle*{10}}
  \put(840,250){\circle*{10}} \put(880,250){\circle*{10}}
  \put(240,780){\circle*{10}} \put(915,250){\circle*{10}}
  \put(340,250){\line( 1, 0){120}} \put(1215,410){\circle*{10}}
  \put(1255,410){\circle*{10}} \put(1290,410){\circle*{10}}
  \put(1215,250){\circle*{10}} \put(1255,250){\circle*{10}}
  \put(1290,250){\circle*{10}} \put(1505,335){\circle*{28}}
  \put(1405,410){\circle*{30}} \put(1405,250){\circle*{30}}
  \put(1405,410){\line( 4,-3){100}} \put(1405,250){\line( 6, 5){ 90}}
  \put(1310,410){\line( 1, 0){100}} \put(1310,250){\line( 1, 0){100}}
  \put(995,410){\circle{28}} \put(1120,410){\circle*{30}}
  \put(930,410){\line( 1, 0){ 55}} \put(1010,410){\line( 1, 0){ 95}}
  \put(1135,410){\line( 1, 0){ 60}} \put(1000,250){\circle{28}}
  \put(1125,250){\circle*{30}} \put(935,250){\line( 1, 0){ 55}}
  \put(1015,250){\line( 1, 0){ 95}} \put(1140,250){\line( 1, 0){ 60}}
  \put(240,410){\circle*{10}} \put(280,410){\circle*{10}}
  \put(315,410){\circle*{10}} \put(600,250){\circle{28}}
  \put(600,410){\circle{28}} \put(105,410){\circle{30}}
  \put(475,410){\circle{30}} \put(710,410){\circle{28}}
  \put(105,250){\circle{30}} \put(475,250){\circle{30}}
  \put(710,250){\circle{28}} \put(490,410){\line( 1, 0){ 95}}
  \put(615,410){\line( 1, 0){ 80}} \put(720,410){\line( 1, 0){100}}
  \put(120,410){\line( 1, 0){100}} \put(340,410){\line( 1, 0){120}}
  \put(490,250){\line( 1, 0){ 95}} \put(615,250){\line( 1, 0){ 80}}
  \put(720,250){\line( 1, 0){100}} \put(120,250){\line( 1, 0){100}}

  \put(105,230){$\underbrace{\hspace*{10.5cm}}_{\hbox{$q_{\nu}$
        vertices}}$} \put(600,130){The $\fR$-index of
    ${^{\gs_{\nu}}G}'(\fR)$}

\end{picture}

\vspace*{1.5cm} From this we see that the boundary component is of type
$\bf I_{\hbf{p$_{\nu}$-d,q$_{\nu}$-d}}$.

There are basically two quite different cases at hand; the first is that
the boundary components are positive-dimensional, the second occurs when
the boundary components reduce to points. The former can be easily handled
with the same methods as above, by splitting off orthogonal complements.
The real interest is in the latter case, and here a basic role is played by
the {\it hyperbolic planes}, which have been dealt with in detail in
\cite{hyp}. We will essentially reduce the rank one case (at least for the
classical groups) to the case of hyperbolic planes, then we explain how the
results of \cite{hyp} apply to the situation here.

\subsection{Positive-dimensional boundary components}
Let $G', P'$ be as above, and consider the hermitian Levi factor
${G'}^{(1)}=M'L'$, which is defined over $k$. Over $\fR$ the factors
$M'(\fR)$ and $L'(\fR)$ are defined. In this section we consider the
situation that the boundary component $F'$ of $\cD'$ defined by $P'$ (i.e.,
$P'(\fR)=N(F')$) is positive-dimensional, or equivalently, that the
hermitian Levi
factor $L'(\fR)$ is non-trivial. As above, we get the following $k$-group
\begin{equation}\label{E23}
N':={G'}^{(1)}\times \cZ_{G'}({G'}^{(1)}).
\end{equation}

The same calculation as in ({\ref{e19.1}) shows that the domain $\cD_{N'}$
defined by $N'$ is the same as that defined by $L'(\fR)\times
\cZ_{G'(\fR)}(L'(\fR))$. Taking the subgroup $N=Res_{k|\fQ}N'$ defines a
subdomain $\cD_N\inn \cD$, which is a product $\cD_N=\cD_{N,\gs_1}\times \cdots
\times \cD_{N,\gs_f}$. Each factor $\cD_{N,\gs}$ is determined by the
corresponding factor of ${^{\gs}L}'(\fR)$. The $\fR$-groups $N'(\fR)$ and
  $N(\fR)$ are determined in terms of the data $D,V,h$ as follows.
\begin{itemize}\item[(I)] If $F'\cong {\bf I_{\hbox{\scsi\bf p-d,q-d}}}$,
  then $\cD_{N'}\cong {\bf I_{\hbox{\scsi\bf p-d,q-d}}}\times {\bf
    I_{\hbox{\scsi\bf d,d}}}$. Note that in terms of the hermitian forms,
  this amounts to the following. Since $h$ has Witt index 1, the maximal
  totally isotropic subspaces are one-dimensional. Let $H_1=<v>$ be such a
  space; we can find a vector $v'\in V$ such that $H=<v,v'>$ is a
  hyperbolic plane, that is, $h_{|H}$ has Witt index 1. It follows that
  $h_{|H^{\perp}}$ is anisotropic. Consider the subgroup
\begin{equation}\label{e22.0} N_k:=\{g\in U(V,h) | g(H)\inn H\}.
\end{equation}
It is clear that for $g\in N_k$, it automatically holds that
 $g(H^{\perp})\inn H^{\perp}$, hence
\begin{equation}\label{e22.1} N_k\cong U(H,h_{|H})\times
  U(H^{\perp},h_{|H^{\perp}}).
\end{equation}
The first factor has $\fR$-points $U(H,h_{|H})(\fR)\cong U(d,d)$, while the
second fulfills $U(H^{\perp},h_{|H^{\perp}})(\fR)\cong U(p-d,q-d)$. Thus
$N_k\cong N'$ as in (\ref{E23}). At any rate, this gives us subdomains of
type
$${\bf I_{\hbox{\scsi\bf d,d}}\times I_{\hbox{\scsi\bf p-d,q-d}}\inn }
\cD_{N'},$$
which, in case $d=p=q$ is the whole domain; in all other cases it is a
genuine subdomain as listed in Table \ref{T1}, defined over $k$, and
$(N',P')$ are incident. It follows from this that $(N,P)$ are incident
over $\fQ$. The components $N_{\gs}(\fR)$ of $N_{\bfb}(\fR)$ are determined
as follows. Let $(p_{\nu},q_{\nu})$ be the signature of $h_{\nu}$ (so that
$p_{\nu}+q_{\nu}=dm$ for all $\nu$). This implies
\[ G(\fR)\cong \prod_{\nu}SU(p_{\nu},q_{\nu}).\]
For each factor, we have the boundary component $F_{\gs}\cong SU(p_{\nu}-d,
q_{\nu}-d)/K$, and for each factor for which $q_{\nu}>d$ this is
positive-dimensional. As above, this leads to subdomains, in each factor,
of type $\bf I_{\hbf{d,d}}\times I_{\hbf{p$_{\nu}$-d,q$_{\nu}$-d}}$, so that in
sum
\begin{equation}\label{E22a} \cD_N\cong \prod_{\nu}\cD_{\nu},\quad
  \hbox{$\cD_{\nu}$ of type $\bf I_{\hbf{d,d}}\times
  I_{\hbf{p$_{\nu}$-d,q$_{\nu}$-d}}$}.
\end{equation}

\item[(II)] Here rank 1 means we have the following $k$-index, $D^{(2)}_{n,1}$
\setlength{\unitlength}{0.004500in}%
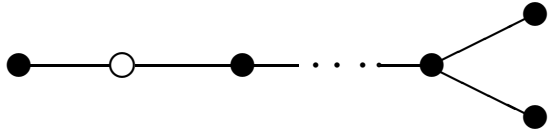
$$\begin{picture}(628,148)(26,666)
\thicklines
\put(300,740){\circle*{28}}
\put(520,740){\circle*{28}}
\put(640,800){\circle*{28}}
\put(640,680){\circle*{28}}
\put(440,740){\circle*{6}}
\put(385,740){\circle*{6}}
\put(410,740){\circle*{6}}
\put( 40,740){\circle*{28}}
\put(160,740){\circle{28}}
\put(525,735){\line( 2,-1){114}}
\put( 45,740){\line( 1, 0){100}}
\put(175,740){\line( 1, 0){115}}
\put(310,740){\line( 1, 0){ 55}}
\put(460,740){\makebox(0.4444,0.6667){\circle{1}}}
\put(460,740){\line( 1, 0){ 55}}
\put(520,740){\line( 2, 1){120}}
\end{picture}$$
In particular, the boundary component is of type $\bf II_{\hbox{\scsi\bf
    n-2}}$ if $\cD'$ is of type $\bf II_{\hbox{\scsi\bf n}}$. This means
also that the ``mixed cases'' only can occur if $\cD'$ is of type $\bf
II_{\hbox{\scsi\bf 4}}$, for then $\bf II_{\hbox{\scsi\bf 2}}\cong$
one-dimensional disc. Of course $\bf II_{\hbox{\scsi\bf
4}}\cong IV_{\hbox{\scsi\bf 6}}$ anyway, so we can
conclude from this that mixed cases do not occur in the hermitian symmetric
setting (for $\fQ$-simple $G$ of rank 1).
The domain $\cD_{N'}$ defined by $N'$ is of
type $\bf II_{\hbox{\scsi\bf n-2}}\times II_{\hbox{\scsi\bf 2}}$. The
components $N_{\gs}(\fR)$ of $N_{\hbf{1}}(\fR)$ are all of type
$U(n-2,\fH)\times U(2,\fH)\inn U(n,\fH)$, so the domain $\cD_N$ is ot
type
\begin{equation}\label{E22a.1} {\bf (II_{\hbf{n-2}}\times II_{\hbf{2}})}^f.
\end{equation}
\item[(III)] Here rank 1 implies the index is one of $C^{(1)}_{1,1}$ (which
  we have excluded) or $C^{(2)}_{n,1}$. The corresponding
  boundary components in these cases are of type $\bf III_{\hbf{n-2}}$.
  The case $C^{(2)}_{2,1}$, for which the boundary component is a point,
  will be
  dealt with later, the others give rise to a subdomain of type $\bf
  III_{\hbox{\scsi\bf 2}}\times III_{\hbox{\scsi\bf n-2}}$. Consequently,
  $\cD_N$ is of type ${\bf (III_{\hbf{n-2}}\times III_{\hbf{2}})}^f$,
  $f=[k:\fQ]$.
\item[(IV)] Here we just have a symmetric bilinear form of Witt index
  1. The $k$-index in this case is necessarily of the form
\setlength{\unitlength}{0.004500in}%
$$\begin{picture}(359,28)(86,766)
\thicklines
\put(240,780){\circle*{28}}
\put(355,780){\circle*{10}}
\put(100,780){\circle{28}}
\put(400,780){\circle*{10}}
\put(245,780){\line( 1, 0){ 85}}
\put(440,780){\circle*{10}}
\put(115,780){\line( 1, 0){115}}
\end{picture}$$
The corresponding boundary component is a point, a case to be
considered below. Splitting off
an anisotropic vector (defined over $k$)
in this case yields a codimension one
subspace $H^{\perp}$ on which $h$ still has Witt index 1, hence the
stabilizer $N'$ defines a subdomain $\cD_{N'}$ of type $\bf
IV_{\hbox{\scsi\bf n-1}}$. $\cD_N$ is then of type ${(\bf
  IV_{\hbf{n-1}})}^f$.
\item[(V)] The only index of rank 1 is
$$
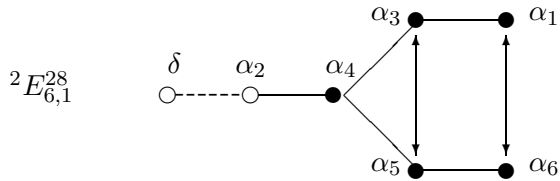
\unitlength1cm
\begin{picture}(14,3)

\put(1,2){$^2E_{6,1}^{28}$}
\put(3,0){
\put(.1,2){\circle{0.2}}
\multiput(.2,1.92)(.15,0){6}{-}
\put(.1,2.3){$\gd$}
\put(1.2,2){\circle{0.2}}
\put(1.3,2){\line(1,0){0.9}}
\put(1,2.3){$\ga_2$}
\put(2.2,2.3){$\ga_4$}

\put(2.3,2){\circle*{0.2}}
\put(2.45,2){\line(1,1){1}}
\put(2.45,2){\line(1,-1){1}}
\put(2.8,3){$\ga_3$}
\put(3.4,3){\circle*{0.2}}
\put(3.4,2){\vector(0,1){.8}}
\put(2.8,1){$\ga_5$}
\put(3.4,1){\circle*{0.2}}
\put(3.4,2){\vector(0,-1){.8}}
\put(3.5,3){\line(1,0){1}}
\put(3.5,1){\line(1,0){1}}
\put(4.6,3){\circle*{0.2}}
\put(4.6,2){\vector(0,1){.8}}
\put(4.9,3){$\ga_1$}
\put(4.6,1){\circle*{0.2}}
\put(4.6,2){\vector(0,-1){.8}}
\put(4.9,1){$\ga_6$}

   }

\end{picture}$$
The vertex denoted $\ga_2$ gives rise to the five-dimensional boundary
component. If $\gd$ denotes the lowest root, then, as is well known, $\gd$
is isotropic (does not map to zero in the $k$-root system), so the root
$\gd$ defines a $k$-subalgebra
$\nn_{\gd}:=\Gg^{\gd}+\Gg^{-\gd}+[\Gg^{\gd},\Gg^{-\gd}]\inn \Gg'$ which is
split over $k$. On the other hand the anisotropic kernel $\cK$ is of type
${^2A}_5$, and $\cK(\fR)\cong U(5,1)$. Clearly $\cK$ and the $k$-subgroup
$N_{\gd}$ defined by $\nn_{\gd}$ are orthogonal, so we get a $k$-subgroup
\[N'=N_{\gd}\times \cK,\]
both factors being defined over $k$. The set of $\fR$-points is then of type
$N'(\fR)\cong SL_2(\fR)\times SU(5,1)$, and the subdomain $\cD_{N'}$ is
\[\cD_{N'}\cong {\bf I_{\hbox{\scsi\bf 1,1}}\times I_{\hbox{\scsi\bf
      5,1}}}.\]
This is one of the domains listed in Table \ref{T1}, incident to the
five-dimensional boundary component. It follows that $\cD_N$ is a product
of factors of this type.
\item[(VI)] There are no indices of hermitian type with rank one for $E_7$.
\end{itemize}
We sum up these results in the following.
\begin{proposition}\label{p23.1} If the rational boundary components for
  $G'$ are positive-dimensional, then Theorem \ref{t12.1} holds for
  $G$. The subdomains defined by the symmetric subgroups $N'\inn G'$ are:
\begin{itemize}\item[(I)] $\bf I_{\hbox{\scsi\bf d,d}}\times
  I_{\hbox{\scsi\bf p-d,q-d}}$.
\item[(II)] $\bf II_{\hbox{\scsi\bf n-2}}\times II_{\hbox{\scsi\bf 2}}$.
\item[(III)] $\bf III_{\hbox{\scsi\bf n-2}}\times III_{\hbox{\scsi\bf 2}}$.
\item[(IV)] $\bf IV_{\hbox{\scsi\bf n-1}}$ (here there are no
  positive-dimensional boundary components).
\item[(V)] $\bf I_{\hbox{\scsi\bf 1,1}}\times I_{\hbox{\scsi\bf 5,1}}$.
\end{itemize}
Note here $\bf I_{\hbox{\scsi\bf 1,1}}\cong II_{\hbox{\scsi\bf 2}} \cong
III_{\hbox{\scsi\bf 1}}\cong IV_{\hbox{\scsi\bf 1}}$.
The corresponding domains $\cD_N$ in $\cD$ defined by the subgroups $N$ are
products of domains of the types listed above.
\end{proposition}

\subsection{Zero-dimensional boundary components}
The restrictions rank equal to one and zero-dimensional boundary componants
are only possible for the domains of type $\bf I_{\hbox{\scsi\bf p,q}},\
III_{\hbox{\scsi\bf 2}}$ and $\bf IV_{\hbox{\scsi\bf n}}$ (see Lemmas
\ref{l5.8.1} and \ref{l5.8.2}). Of these, the
last case requires no further discussion: as above we find a codimension
one $k$-subspace $V'\inn V$, on which $h$ still is isotropic, and take its
stabilizer as $N'$. This gives a $k$-subgroup $N'\inn G'$, and defines a
subdomain $\cD_{N'}$ of type $\bf IV_{\hbox{\scsi\bf n-1}}$. In the
${^2A}^{(d)}$ case we may assume $d\geq 3$: the $d=1$ case is again easily
dealt with as above. We have a $K$-vector space $V$ ($K|k$ imaginary
quadratic) of dimension $p+q$ and a ($K$-valued) hermitian form $h$ of Witt
index 1 on $V$. By taking a $K$-subspace $V'\inn V$ of codimension one,
such that $h_{|V'}$ still has Witt index 1, we get the $k$-subgroup $N'$ as
the stabilizer of $V'$. Then the domain $\cD_{N'}$ is either of type $\bf
I_{\hbox{\scsi\bf p-1,q}}$ or $\bf I_{\hbox{\scsi\bf p,q-1}}$, and by
judicious choice of $V'$ we can assume the first case, which is the domain
listed in Table \ref{T2}. The $d=2$ case is ``lifted'' from the
corresponding $d=2$ case with involution of the first kind: if $D$ is
central simple of degree 2 over $K$ with a $K|k$-involution, then
(\cite{A}, Thm.~10.21) $D=D'\otimes_kK$, where $D'$ is central simple of
degree 2 over $k$ with the canonical involution. Consequently,
$$U(V,h)=U(V'\otimes_kK,h'\otimes_kK)=U(V',h')_K,$$
the group is just the group $U(V',h')$ lifted to $K$. Since $U(V',h')$ has
index $C^{(2)}_{n,1}$, while $U(V',h')_K$ has index $A^{(2)}_{2n-1,1}$, it
follows that the boundary component is a point only if $n\leq 2$. This
implies that if $d=2$, the index is $A^{(2)}_{3,1}$, the domain is $\bf
I_{\hbf{2,2}}\cong IV_{\hbf{4}}$, so $U(V',h')_K$ is isomorphic to an
orthogonal group over $k$ in six variables. As we just saw, in this case
there is a subdomain defined over $k$ of type $\bf IV_{\hbf3}\inn
IV_{\hbf4}$. So we assume $d\geq 3$. Then, as we have seen, the
boundary component
$F'\cong \bf I_{\hbox{\scsi\bf p-d,q-d}}$ will be
zero-dimensional $\iff$ $q=d$ (respectively $F\cong {\bf
  I_{\hbf{p$_1$-d,q$_1$-d}}\times \cdots \times I_{\hbf{p$_f$-d,q$_f$-d}}}$
   will be
  zero-dimensional $\iff$ $q_{\nu}=d,\ \forall_{\nu}$.
Here there are two possibilities:
\begin{itemize}\item[1)] $p=q=d$, the group $N_k$ of (\ref{e22.0}) is
  $N_k\cong G'$. This is the case of {\it hyperbolic planes}.
\item[2)] $p>q=d$, the group $N_k$ of (\ref{e22.0}) is over $\fR$ just
  $N_k(\fR) = U(d,d)\times U(p-d)\inn U(p,d)\cong G'(\fR)$.
\end{itemize}
Note that in the second case the domain $\cD_{N_k}$ defined by $N_k$ is of
type $\bf I_{\hbox{\scsi\bf d,d}}$, a maximal tube domain in $\bf
I_{\hbox{\scsi\bf p,q}}$. So we are also finished in this case. For
completeness, let us quickly go through the details to make sure
nothing unexpected happens.
\begin{proposition}\label{p24.1} Let $G'$ have index ${^2A}^{(d)}_{n,1},\
  d=q,\ p>q$, $n+1=p+q$,
  and let $P'$ denote the corresponding standard parabolic and
  $N'= N_k$, where $N_k\inn G'$ the symmetric subgroup defined in
  (\ref{e22.0}), where $H$
  is the hyperbolic plane spanned by the vector which is stabilized by $P'$
  and its ortho-complement ($v'$: $h(v,v')=1$). Then $(P',N')$ are
  incident, in fact standard incident.
  Consequently, $P=Res_{k|\fQ}P'$ and $N=Res_{k|\fQ}N'$ are
  incident over $\fQ$.
\end{proposition}
{\bf Proof:} We know that $N'(\fR)\cong U(q,q)\times U(p-q)$
which gives rise to the
maximal tube subdomain $\bf I_{\hbox{\scsi\bf q,q}}\inn I_{\hbox{\scsi\bf
    p,q}}$ of Table \ref{T2}. We need to check that the standard boundary
component $F'$ stabilized by $P'(\fR)$ is also a standard boundary
component of $\cD_{N'}$; in particular we need the common maximal
$\fR$-split torus in $P'$ and $N'$. This is seen in (\ref{e22.0}), the
$\fR$-split torus being contained in the hermitian Levi factor of
$P'(\fR)$, which is contained in $N'(\fR)$.
Consider the group $P'\cap N'$; this
is nothing but the stabilizer of $v$ in $H$, which is a maximal standard
parabolic in $N'$. Since $\nu$ determines the boundary component $F$, both
in $G'(\fR)$ and in $N'(\fR)$, it is clear that $F$ is a boundary component
of $\cD_{N'}$. It follows that $(P',N')$
are incident, and this implies (see the discussion preceeding Proposition
\ref{p19.2}) that $(P,N)\inn G$ are incident. \ende
We are left with the following cases: $\bf III_{\hbox{\scsi\bf 2}}$ with
index $C^{(2)}_{2,1}$ and $\bf I_{\hbox{\scsi\bf q,q}}$ with index
${^2A}^{(d)}_{2d-1,1},\ d\geq3$.
These indices are described in terms of hermitian
forms as follows. Let $D$ be a central simple division algebra
over $K$ ($K=k$ for $d=2$ and $K|k$ is imaginary quadratic if $d\geq3
$) and assume further that $D$ has a $K|k$-involution, $V$ is a
two-dimensional right vector space over $D$ and $h:V\times V\lra D$ is a
hermitian form which is isotropic. Then $d=2$ gives groups with index
$C^{(2)}_{2,1}$, and $d\geq 3$ gives groups with indices
${^2A}^{(d)}_{2d-1,1}$.
\begin{lemma}\label{l25.1} There exists a basis $v_1,v_2$ of $V$ over $D$
  such that the form $h$ is given by $h({\bf x},{\bf
    y})=x_1\overline{y}_2+x_2\overline{y}_1,\ {\bf x}=(x_1,x_2),\ {\bf
    y}=(y_1,y_2)$.
\end{lemma}
{\bf Proof:} Let $v$ be an isotropic vector, defined over $k$. Then there
exists an isotropic
 vector $v'$, such that $h(v,v')=1$, hence also
$h(v',v)=1$. Let ${v'}=({v'}_1,{v'}_2)$, and set
$\gd={v'}_1\overline{v'}_2$, so that
$h(v',v')=\gd+\overline{\gd}$. Then the matrix of $h$ with respect to
the basis $v,v'$ is $H'={0\ 1\choose 1\ \ge}$, where
$\ge=\gd+\overline{\gd}$. Now setting
$$w=(w_1,w_2)=(-v_1\overline{\gd}+v_1',-v_2\overline{\gd}+v_2')$$
we can easily verify $h(w,w)=0,\ \
h(v,w)=h(w,v)=1$. Since the change of basis transformation is defined over
$k$, the matrix of the hermitian form with respect
to this $k$-basis $v,w$ is $H={0\ 1\choose 1\ 0}$. \ende
So as far as the $\fQ$-groups are concerned, we may take the standard
hyperbolic form given by the matrix $H$ as defining the hermitian form on
$V$. We remark that the situation changes when one considers arithmetic
groups, but that need not concern us here. At any rate, a two-dimensional
right $D$-vector space $V$ with a hermitian form as in Lemma \ref{l25.1} is
what we call a {\it hyperbolic plane}, and this case was studied in detail
in \cite{hyp}. There it was determined exactly what kind of
symmetric subgroups
exist. These derive from the existence of splitting subfields $L\inn D$,
which may be taken to be cyclic of degree $d$ over $K$, if $D$ is central
simple of degree $d$ over $K$. In fact, we have subgroups (\cite{hyp},
Proposition 2.4) $U(L^2,h)\inn U(D^2,h)$, which give rise to the following
subdomains:
\begin{itemize}
\item[1)] $d=2$; $\cD_L\cong \left(\begin{array}{cc}\tau_1 & 0 \\ 0 &
      b^{\zeta_1}\tau_1
\end{array}\right)\times \cdots \times \left(\begin{array}{cc}\tau_1 & 0 \\ 0 &
  b^{\zeta_f}\tau_1
\end{array}\right)$, where $\zeta_i:k\lra {\Bbb R}$ denote the distinct real
embeddings of $k$.
\item[2)] $d\geq 3$; $\cD_L\cong \left(\begin{array}{ccc}\tau_1 & & 0 \\ &
      \ddots & \\ 0 & & \tau_d\end{array}\right)^f$.
\end{itemize}
In other words, for hyperbolic planes we find subdomains of the following
kinds
\begin{equation}\label{E25}
\bf III_{\hbox{\scsi\bf 1}}\inn III_{\hbox{\scsi\bf 2}},\quad
(I_{\hbox{\scsi\bf 1,1}})^d\inn I_{\hbox{\scsi\bf d,d}}.
\end{equation}
The latter one is a polydisc, coming from a maximal set of strongly orthogonal
roots, i.e., satisfying 1), 2'') and 3'').
The first case is the only exception to the rule that we have
symmetric subgroups $N'\inn G'$ with $\rank_{\fR}N'=\rank_{\fR}G'$.

\vspace*{.2cm}
\noindent{\bf Proof of Theorem \ref{t12.1}:} We have split the set of cases
up into the three considered in \S4, 5 and 6. Corollary \ref{c14.1} proves
\ref{t12.1} for the split over $\fR$ case and Proposition \ref{p19.2} for
the rank $\geq 2$ case and positive-dimensional boundary components.
For $\rank \geq 2$ and zero-dimensional boundary components, Proposition
\ref{p23} shows that with the exception given Theorem \ref{t12.1}
holds in this case also. In the case of rank 1, Proposition \ref{p23.1}
verifies \ref{t12.1} for the case that the boundary components are
positive-dimensional, and Proposition \ref{p24.1} took care of the rest of
the cases excepting hyperbolic planes. Then the results of \cite{hyp}
verify \ref{t12.1} for hyperbolic planes, thus completing the proof. \ende

\vspace*{.2cm}
\noindent{\bf Proof of the Main Theorem:} The first statement is covered by
Theorem \ref{t12.1}. The statements on the domains for the exceptions
follow from (\ref{E20}) and (\ref{E25}). It remains to consider the
condition 4). This is fulfilled for the groups $N$ utilized above
by construction. For the
exceptional cases this is immediate, as we took subgroups defined by
symmetric closed sets of roots. Let us sketch this again for the classical
cases, utilizing the description in terms of $\pm$symmetric/hermitian
forms. The objects $D,\ V,\ h$ and $G'$ will have the meanings as
above. Let $s=\rank_{k}G'$, and let $H_s$ be an $s$-dimensional (maximal)
totally isotropic subspace in $V$, with basis $h_1,\ldots, h_s$. Let
$h_i'\in V$ be vectors of $V$ with $h(h_i,h_j')=\gd_{ij}$,
$H_s'=<h_1',\ldots, h_s'>$ and set $H=H_s\oplus H_s'$. Then $h_{|H}$ is
non-degenerate of index $s$, and $H$ splits into a direct sum of hyperbolic
planes, $H=V_1\oplus \cdots \oplus V_s$. The form $h$ restricted to
$H^{\perp}$ is anisotropic; the semisimple anisotropic kernel is
$SU(H^{\perp},h_{|H^{\perp}})$. Fixing the basis $h_1,\ldots,
h_s,h_1',\ldots, h_s'$ for $H$ amounts to the choice of maximal $k$-split
torus $S'$. For each real prime $\nu$, $(H_{\nu},h_{\nu})$ is a
$2ds$-dimensional $\fR$-vecotr space with $\pm$symmetric/hermitian
form. Choosing an $\fR$-basis of $H_{\nu}$ amounts to choosing a maximal
$\fR$-split torus of $SU(H_{\nu},h_{\nu})$, and a choice of basis for a
maximal set of hyperbolic planes (over $\fR$) amouts to the choice of
maximal $\fR$-split torus. Similarly,
$(H_{\nu}^{\perp},{h_{\nu}}_{|H_{\nu}^{\perp}})$ is an $\fR$-vector space,
$h_{|{H^{\perp}_{\nu}}}$ has some index $q_{\nu}$,
and one can find a maximal set of
hyperbolic planes $W_1,\ldots,W_r$, such that
$H_{\nu}^{\perp}=(W_1)_{\nu}\oplus \cdots \oplus (W_r)_{\nu}\oplus W'$,
where ${h_{\nu}}_{|W'}$ is anisotropic over $\fR$. A choice of basis of the
$(W_i)_{\nu}$ amounts to the choice of maximal $\fR$-split torus, and a
choice of basis, over $\fR$, of $V_{\nu}$ amounts to the choice of maximal
torus defined over $\fR$. From these descriptions we see that the polydisc
group $N_{\Psi}$ defined by the maximal set of strongly orthogonal roots
$\Psi$ splits into
a component in $SU(H,h_{|H})$ and a component in
$SU(H^{\perp},h_{|H^{\perp}})$, say $N_{\Psi}=N_{\Psi,1}\times
N_{\Psi,2}$. Then $N_{\Psi,2}\inn SU(H^{\perp},h_{|H^{\perp}})(\fR)$ and
$N_{\Psi,1}\inn SU(H,h_{|H})$. Since the subgroup
$SU(H^{\perp},h_{|H^{\perp}})$ is contained in all the groups $N$ we have
defined, we need only consider $N_{\Psi,1}$. $H$ is a direct sum of
hyperbolic planes $V_i$, and the question is whether the corresponding
polydisc group is contained in $SU(V_i,h_{|V_i})$. But this is what was
studied in \cite{hyp}; the answer is affirmative. It follows that with the
one exception stated, $C^{(2)}_{2,1}$, $N_{\Psi}\inn N$. \ende

\end{document}